\documentclass[11pt]{article}
\usepackage{microtype}
\usepackage{enumitem}
\usepackage[a4paper,margin=2.5cm]{geometry}
\usepackage{amsmath,amssymb,amsthm,mathtools}
\usepackage{bm}
\usepackage{xcolor}
\usepackage{authblk}
\usepackage{booktabs}
\usepackage{tabularx}
\usepackage{array}
\usepackage{orcidlink}
\usepackage{graphicx}
\usepackage{subcaption}
\usepackage{comment}
\usepackage{makecell}
\usepackage{hyperref}
\usepackage[section]{placeins}
\usepackage{booktabs}

\usepackage{algorithm}
\usepackage{algpseudocode}
\usepackage{float}

\definecolor{citationblue}{RGB}{0,70,140}

\title{Adaptive Local Proper Orthogonal Decomposition via Graph Partitioning and Provable Error-Bounded Lossy Compression for Parametric Series in Finite-Dimensional Function Spaces}

\author[1]{Roberto Nuca\,\orcidlink{0000-0002-9031-5668}%
  \thanks{Corresponding author: roberto.nuca@kaust.edu.sa}}
\author[1,2]{Matteo Parsani\,\orcidlink{0000-0001-7300-1280}}
\author[1]{Lisandro Dalcin\,\orcidlink{0000-0001-8086-0155}}

\affil[1]{Computer, Electrical and Mathematical Sciences and Engineering Division,
King Abdullah University of Science and Technology,
Thuwal, 23955-6900, Saudi Arabia}

\affil[2]{Physical Science and Engineering Division, King Abdullah University of Science and Technology,
Thuwal, 23955-6900, Saudi Arabia}

\date{\today}

\newtheorem{theorem}{Theorem}[section]

\newtheorem{proposition}[theorem]{Proposition}

\newtheorem{corollary}[theorem]{Corollary}
\theoremstyle{remark}
\newtheorem{remark}{Remark}[section]

\newcommand{\vstack}{\operatorname*{vstack}}
\newcommand{\rank}{\operatorname{rank}}

\newcommand{\NRMSE}{\mathrm{NRMSE}}
\newcommand{\CR}{\mathrm{CR}}

\begin{document}
\maketitle

\begin{abstract}
We introduce a spatially adaptive proper orthogonal decomposition (POD) method for error-bounded lossy compression of snapshot data in finite-dimensional function spaces. The method combines hierarchical graph partitioning with independent low-rank approximations on spatial subdomains, allowing the local approximation ranks to adapt to the spatial distribution of information in the data. For partitions that are orthogonal in the discrete inner product, we establish an exact decomposition of the global squared weighted reconstruction error into local contributions and characterize the weighted-norm optimality and coefficient-space quasi-optimality of local POD approximations for fixed partitions and ranks. These results provide a rigorous mechanism for allocating a prescribed global error tolerance among subdomains and selecting the smallest admissible local ranks. We further develop a sufficient error bound for the coefficient-space normalized root-mean-square error (NRMSE). A greedy hierarchical refinement algorithm accepts a spatial subdivision only when it reduces the storage required by the local low-rank factors while preserving the prescribed reconstruction-error bound. The method is assessed using synthetic data and snapshot sequences from two- and three-dimensional compressible-flow simulations. The numerical results illustrate how spatially adaptive ranks and partitioning can improve compression relative to global POD and fixed spatial partitions, particularly when the data exhibit spatially heterogeneous dynamics.
\end{abstract}



\section{Introduction}
\label{sec:introduction}
Large-scale numerical simulations often produce more data than can be stored, transferred, or analyzed conveniently. After spatial discretization, a time-dependent calculation yields snapshots $u_h(t_1),\ldots,u_h(t_n)\in V_h$, where $V_h$ is a finite-dimensional function space. If $m=\dim(V_h)$, storing the snapshot coefficients alone requires $m\,n$ real numbers, in addition to the mesh and other metadata. For large simulations, this data volume affects checkpointing, visualization, post-processing, data transfer, uncertainty quantification, and surrogate model construction. This motivates lossy compression methods that reduce storage while keeping the reconstruction error below a prescribed tolerance.

Error-bounded scientific compressors exploit different forms of structure. Prediction-based methods such as SZ control pointwise errors; ZFP uses block transforms and supports fixed-rate, fixed-precision, and fixed-accuracy compression; MGARD relies on a multilevel decomposition with computable norm bounds; and TTHRESH exploits multilinear structure through a higher-order singular value decomposition \cite{DiCappello2016SZ,Lindstrom2014ZFP,DiffenderferEtAl2019ZFP,AinsworthEtAl2019MGARD,BallesterRipollEtAl2020TTHRESH}. A recent approach combines randomized singular value decomposition with classical lossy compression algorithm, like SZ, to combine the benefits of both approaches \cite{NucaParsaniTurkiyyah2025}. Here we consider a more specific problem: compressing a sequence of functions represented in the same finite-dimensional space. This setting exploits temporal redundancy while measuring reconstruction error in a norm inherited from the underlying function space.

Singular value decomposition (SVD) and proper orthogonal decomposition (POD) provide a natural basis for this approach. Once a basis of $V_h$ is chosen, the snapshot coefficients form a matrix $\bm{A}=[\boldsymbol{u}_1,\ldots,\boldsymbol{u}_n] \in\mathbb{R}^{m\times n}$. The Eckart--Young--Mirsky theorem states that the rank-$p$ truncated SVD is a best rank-$p$ approximation for unitarily invariant norms and, for the spectral and Frobenius norms, its error is determined by the discarded singular values \cite{EckartYoung1936,Mirsky1960}. In reduced-order modeling, the same construction is commonly interpreted as POD \cite{Sirovich1987,BerkoozHolmesLumley1993}. If the inner product on $V_h$ is represented by a symmetric positive-definite Gram matrix $\bm{M}$, the SVD of $\bm{M}^{1/2}\bm{A}$ yields the POD that is optimal in the corresponding discrete norm.

A single global POD basis can nevertheless be inefficient for spatially localized or multiscale phenomena. A moving front, shock, vortex, or intermittent structure may occupy only a small part of the domain but still require many global modes. Several approaches address the cost of constructing large POD approximations. Domain-decomposed and partitioned methods construct or approximate a global POD basis from distributed data \cite{BeattieEtAl2006,WangMcBeeIliescu2016}; incremental methods update a low-rank representation as snapshots arrive \cite{Brand2006,FareedEtAl2018,FareedSingler2020}; HAPOD organizes local POD computations on a tree and controls the error of the resulting global basis \cite{HimpeLeibnerRave2018}; and randomized methods approximate the dominant range through sketches \cite{HalkoMartinssonTropp2011}. Although these algorithms differ in purpose and guarantees, their partition, stream, or hierarchy is generally used to construct one global reduced space. In our method, by contrast, the spatial partition and the independent local bases are retained in the compressed representation.

A separate line of work retains local reduced spaces. Domain-decomposed reduced-order models have been developed for moving shocks, nonlinear mechanics, and other problems with strongly localized behavior \cite{LuciaKingBeran2002,BaigesCodinaIdelsohn2013,KerfridenEtAl2013}. The domain-decomposition non-intrusive reduced-order model (DDNIROM) of Xiao et al.\ is particularly relevant \cite{XiaoEtAl2019DDNIROM}. DDNIROM constructs local POD bases on prescribed subdomains and learns radial-basis-function maps for the evolution of the local coefficients, using information from neighboring subdomains to represent coupling. Its numerical examples show that local singular values may decay more rapidly than the global ones and that different regions may benefit from different basis dimensions. Its goal, however, is non-intrusive time prediction rather than compression. The partition and local dimensions are selected through physical considerations and retained-energy criteria, rather than through a storage-aware partition search under a certified global reconstruction tolerance.

Spatially adaptive low-rank representations have also been considered in image and tensor compression. Souza and Pedrini combine image partitioning with local SVD approximations and spatially varying ranks \cite{SouzaPedrini2019}. Their method reflects the same basic observation that different regions can have different spectral complexity, but its partitioning criteria and error measures are tailored to image data.

A closer structural connection is provided by the adaptive hierarchical subtensor method of Ehrlacher et al. \cite{EhrlacherEtAl2021AdaptiveTensor}. That method partitions a tensor into subtensors and assigns a separate low-rank tensor approximation to each one. Neither the partition nor the local ranks are prescribed in advance; they are adapted to meet an accuracy requirement while reducing storage. Our problem differs in that the data are snapshots of functions on a spatial discretization, the partition is tied to the physical domain, and the local approximations are local POD factorizations whose errors are measured in the discrete function-space norm.

The recent work of Shamrai is particularly close to the matrix-compression aspect of the present study \cite{Shamrai2026ConcatenatedSVD}. Given a collection of matrices with a common row dimension, Shamrai forms clusters, horizontally concatenates the matrices within each cluster, and applies a shared truncated SVD. Candidate merges are assessed through Weyl-type information, incremental residuals, or an incremental SVD approximation of the Frobenius reconstruction error. There is a direct algebraic connection with a row-partitioned snapshot matrix: after transposition, a parent matrix is the horizontal concatenation of its transposed child blocks. Compressing the parent imposes one temporal subspace on all of those blocks, whereas compressing them separately allows different temporal subspaces in different spatial regions. The two formulations nonetheless remain distinct. Shamrai clusters a collection of pre-existing matrices using a prescribed target rank, whereas we construct a geometrically admissible hierarchy for one spatially discretized snapshot matrix and select the smallest admissible rank separately in each final subdomain.

To our knowledge, the cited approaches do not combine the particular features considered here: a spatially meaningful decomposition of one function-valued snapshot matrix, independent local POD approximations on the final subdomains, local ranks selected from a single global error budget, and a partition criterion based on the storage required by the compressed factors. We generate candidate subdomains through a hierarchy of spatial partitions and retain a separate local POD representation on each accepted leaf. The rank of each local approximation is not fixed in advance but is chosen as the smallest rank compatible with its assigned portion of the global tolerance.

The error analysis assumes a diagonal Gram matrix, which makes every index partition orthogonal in the discrete inner product. Thus, the global squared weighted error is obtained exactly by summing the discarded spectral energies of the local approximations. A global squared-error budget can therefore be divided among the subdomains, and each local truncation rank is selected as the smallest rank satisfying its local budget. When accuracy is reported in the ordinary coefficient-space NRMSE, we use the norm-equivalence bound derived in Section~\ref{sec:localsurrogates}; this coefficient-space certificate is sufficient but generally conservative.

The same hierarchy determines where localization is worthwhile. For each candidate split, the children inherit the error budget of their parent and are compressed at their smallest admissible ranks. The split is retained only when the combined child representation requires less algebraic storage than the parent representation. This comparison accounts for the central trade-off: spatial localization may accelerate singular-value decay and reduce the local ranks, but every additional subdomain introduces another set of temporal factors and partition metadata. For a fixed partition and prescribed local ranks, the local POD approximations are optimal in the weighted Frobenius norm within the corresponding block-separable approximation class and quasi-optimal in the ordinary Frobenius norm. The main contributions are:
\begin{enumerate}
    \item an exact coefficient-space error identity for partitioned SVD, together with an exact weighted-error identity, a coefficient-space error bound, and fixed-partition weighted optimality and Frobenius quasi-optimality for partitioned POD;

    \item a deterministic allocation of a global error budget and a rule that selects the smallest admissible POD rank in each subdomain, supplemented by a sufficient bound for the coefficient-space NRMSE;

    \item a graph-based hierarchical partitioning procedure that retains localization only when it preserves the error certificate and reduces the algebraic storage cost; and

    \item a numerical assessment against classical POD and adaptive local POD on fixed spatial partitions, including the effects of local spectral decay, error allocation, and partition refinement on compression.
\end{enumerate}

The reconstruction-error guarantee is global. The partition search is restricted to the hierarchy generated by the graph-partitioning procedure and, in its present top-down greedy form, is not claimed to minimize storage over all possible spatial partitions.

The remainder of the paper is organized as follows.
Section~\ref{sec:function-space} introduces the finite-dimensional function-space setting and compatible spatial partitions.
Section~\ref{sec:localsurrogates} develops the local surrogates, the error identities, and the error-budget allocation.
Section~\ref{sec:compression} presents the algebraic storage model and the adaptive partitioning algorithm. 
Sections~\ref{sec:numerical-results} and~\ref{sec:conclusions} report the numerical results and conclusions, respectively.

\paragraph{Notation} Vector and matrix indices are numbered from one. For a matrix $\bm{A}$, the notation $\bm{A}[I_r,I_c]$ denotes the submatrix formed by extracting row indexes specified by the index set $I_r$, and column indexes specified by the index set $I_c$. The symbol ``$:$'' is used to denote all rows or all columns. For example, $\bm{A}[I,:]$ denotes the submatrix obtained from $\bm{A}$ by selecting the row indices specified in $I$ and all columns. When a number is appended to ``$:$'', e.g $\bm{A}[:p,:]$, the index set for rows is ${1,2,\dots,p}$. When a number is prepended to ``$:$'', e.g $\bm{A}[p:,:]$, the index set for rows is ${p+1,\dots,m}$. The same mechanism works for columns. Here we presented an example on matrices; however, the same notation applies to multidimensional arrays, including vectors, matrices, and tensors.

\section{Mathematical Setting}
\label{sec:function-space}
We consider snapshot data represented in a fixed finite-dimensional function space over a fixed computational mesh. We assume the basis functions and the Gram matrix are independent of the snapshot index. This setting covers time-dependent and parametrized partial differential equations on static meshes. A deforming-mesh problem can be included in this framework only after mapping all snapshots to a common reference representation with a fixed basis and Gram matrix. This section introduces the mesh partition, the induced spatial and algebraic partitions, and the finite-dimensional notation used in the error analysis.

\subsection{Mesh Description, Mesh Partitioning and Degrees of Freedom}
\label{sec:meshdescr}
Let $\Omega\subset\mathbb R^d$ be an open set on which the functions of interest are defined. For time-dependent problems, we consider the space-time cylinder $\Omega\times[0,T]$. A \textbf{\textit{mesh}}, or \textbf{\textit{tessellation}}, of $\Omega$ is a finite collection 
\begin{equation*}
    \mathcal T_h:=\{T_j\}_{j=1}^{N_c}
\end{equation*}
of subsets of $\mathbb{R}^d$, called \textbf{\textit{elements}} or \textbf{\textit{cells}}. The symbol $N_c$ denotes the number of elements in the mesh. We require the closures of the elements to coincide with $\overline{\Omega}$ and the cell interiors to be pairwise disjoint: 
\begin{equation*}
    \overline{\Omega} = \bigcup_{j=1}^{N_c}\overline{T_j} \qquad\qquad T_i^\circ\cap T_j^\circ = \varnothing\quad\text{for}\quad i\neq j.
\end{equation*}
Here, 
\begin{equation*}
    \overline{\Omega} = \Omega\cup\partial\Omega
\end{equation*}
indicates the closure of $\Omega$. The closures are used because the mesh must represent both the external boundary $\partial\Omega$ and the interfaces shared by neighboring cells. The shapes of these elements are typically triangles, quadrilaterals, polygons, etc. in two dimensions (2D), and tetrahedra, hexahedra, or general polyhedra in three dimensions (3D). We partition the cells of $\mathcal T_h$ into $N_p$ nonempty, connected groups. To identify the cells in each group, we define
\begin{equation}
    \mathcal C := \{\mathcal C_k\}_{k=1}^{N_p},
    \label{eq:cell-partition}
\end{equation}
where $\mathcal C_k\subset\{1,\ldots,N_c\}$ contains the indices of the cells belonging to the $k$-th group. We call $\mathcal C$ a
\textbf{\textit{cell partition}} if
\begin{equation}
    \varnothing\ne\mathcal C_k\subset\{1,\ldots,N_c\},
    \qquad
    \mathcal C_i\cap\mathcal C_j = \varnothing
    \quad\text{for }i\ne j,
    \qquad
    \bigcup_{k=1}^{N_p}\mathcal C_k = \{1,\ldots,N_c\}.
    \label{eq:cell-partition-properties}
\end{equation}
Furthermore, we require the cells indexed by each $\mathcal C_k$ to form a connected set in the cell-adjacency graph, where two cells are adjacent when they share a mesh interface. We denote the number of cells in the $k$-th group by
\begin{equation*}
    N_k := |\mathcal C_k|.
\end{equation*}
Since every subdomain contains at least one cell, we have
\begin{equation*}
    1\leq N_p\leq N_c,
\end{equation*}
where $N_p=1$ corresponds to one group of cells, and $N_p=N_c$ means that each cell is a subdomain and each subdomain is a single cell. The $k$-th \textbf{\textit{subdomain}} is the open set
\begin{equation}
    \Omega_k := \operatorname{int}\left(
    \bigcup_{j\in\mathcal C_k}\overline{T_j}
    \right).
    \label{eq:subdomain-definition}
\end{equation}
The collection $\{\Omega_k\}_{k=1}^{N_p}$ is called the \textbf{\textit{spatial partition}} induced by the cell partition $\mathcal C$. The subdomains satisfy
\begin{align}
    \overline{\Omega}
    &=
    \bigcup_{k=1}^{N_p}\overline{\Omega_k},
    \label{eq:prequirements}\\
    \Omega_i\cap\Omega_j
    &=
    \varnothing
    \quad\text{for }i\ne j.
    \label{eq:prequirements2}
\end{align}
The closures appear in Eq.~\eqref{eq:prequirements} because the open subdomains do not contain their shared interfaces or the external boundary. Therefore, their closures may intersect along interfaces, whereas the open subdomains themselves are pairwise disjoint.

The cell partition $\mathcal{C}$ can be generated, for example, using a graph-partitioning algorithm. A popular option is \textsc{Metis}~\cite{karypis1998metis}, which has been successfully adopted for partitioning unstructured meshes in the context of finite elements. We remark that the subsequent analysis is independent of the particular mesh partitioning algorithm, provided that the cell groups are nonempty, pairwise disjoint, connected, and collectively exhaustive. We also remark that the graph partitioning algorithm can account for edge/vertex weights. This is particularly relevant when information per cell is not uniform over the mesh. Indeed, when data are represented over the physical domain, for example, via discretization of space-time varying functions, the total number of degrees of freedom might not coincide with $N_c$. Thus, we next associate the spatial partition with a partition of the algebraic degrees of freedom (DoFs). Let $m$ denote the total number of scalar DoFs in one discrete state. If the state contains several physical variables, each scalar coefficient is counted separately in $m$. A DoF is a scalar coefficient and does not necessarily correspond to a physical point. For example, a modal coefficient is associated with an element rather than with a collocation point. The relation between $m$ and $N_c$ depends on the discretization approach. In particular, $m$ can be equal to, or larger than $N_c$. The case $m=N_c$ corresponds to one scalar value per cell. Therefore, no general ordering between these quantities is assumed. We identify the DoFs with the index set $\{1,\ldots,m\}$, and we group the degree-of-freedom indices according to the subdomains and define
the corresponding \textbf{\textit{index partition}} by
\begin{equation*}
    \mathcal P := \{I_k\}_{k=1}^{N_p},
\end{equation*}
where the nonempty index sets satisfy
\begin{equation*}
    \varnothing\ne I_k\subset\{1,\ldots,m\},
    \qquad
    I_i\cap I_j = \varnothing
    \quad\text{for }i\ne j,
    \qquad
    \bigcup_{k=1}^{N_p}I_k = \{1,\ldots,m\}.
\end{equation*}

\begin{figure}[H]
    \centering
    \includegraphics[width=0.65\linewidth]{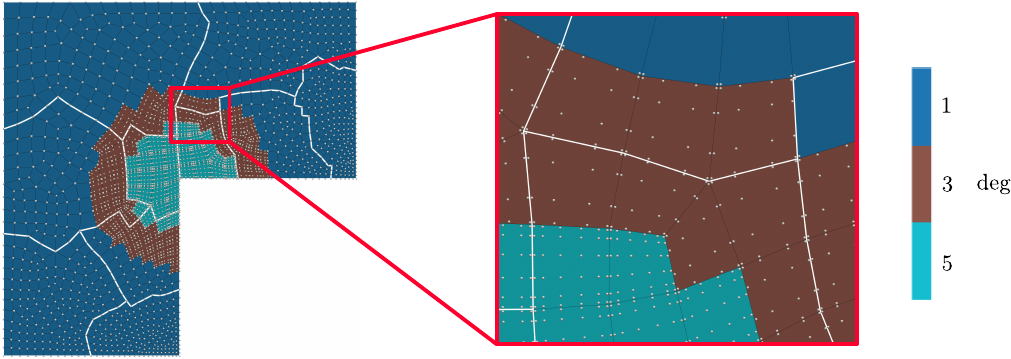}
    \caption{Example of weighted mesh partitioning of an unstructured mesh of rectilinear quadrilaterals equipped with discontinuous collocated polynomials of degrees 1, 3, and 5. The mesh is partitioned into ten subdomains $\Omega_1,\dots,\Omega_{10}$; the white lines represent the subdomains' boundaries. Different colors represent different polynomial degrees: 1, 3, and 5 in the example. Gray dots represent the collocation points of each element. The dots at the boundaries of each element are slightly shifted to make clear that multiple degrees of freedom reside at the same physical point. In this case, we assign weights to graph vertices based on the number of DoFs per cell.}
    \label{fig:weighted-partitioning}
\end{figure}

We call $I_k$ an \textbf{\textit{index set}}. For discretizations with element-local DoFs, $I_k$ is obtained by collecting the DoF indices associated with the cells whose indices belong to $\mathcal C_k$. For conforming discretizations with degrees of freedom shared by neighboring cells, a fixed ownership rule is required to keep the sets $I_k$ disjoint. Once this rule has been fixed, each subdomain is associated with one index set. Figure~\ref{fig:weighted-partitioning} shows an example of weighted partitioning. In particular, we created a mesh equipped with a space of discontinuous polynomials of three values of polynomial degree; the figure illustrates how DoFs are distributed over the mesh. The partitioning is weighted by the number of DoFs per cell. Indeed, the number of cells per subdomain varies significantly, but the number of DoFs per subdomain is balanced.

Thus, summarizing, the construction follows the sequence
\begin{equation*}
    \mathcal C_k
    \longrightarrow
    \Omega_k
    \longrightarrow
    I_k,
\end{equation*}
from a group of mesh cells, to a physical subdomain, and finally to a block of algebraic indices. The index partition $\mathcal P$ determines the row blocks of the snapshot matrix introduced below.

\subsection{Spaces of Functions} 
We use the mesh formalism described in the previous subsection to distribute quantitative information across physical space. Let $V$ be a Hilbert with inner product $(\cdot,\cdot)_V$. For example, $V=L^2(\Omega)$ equipped with the inner product
\begin{equation*}
    (u,v)_{L^2(\Omega)} :=
    \int_\Omega u(x)\,v(x)\,\mathrm{d}x,
\end{equation*}
whereas $V=H^1(\Omega)$ equipped with
\begin{equation*}
    (u,v)_{H^1(\Omega)} :=
    \int_\Omega u(x)\,v(x)\,\mathrm{d}x
    +
    \int_\Omega\nabla u(x)\cdot\nabla v(x)\,\mathrm{d}x.
\end{equation*}
For vector-valued states, $V$ can be taken as an appropriate product Hilbert space.

To approximate the functions of $V$, we introduce $V_h$, a finite-dimensional vector subspace of $V$ containing functions that approximate those in $V$. The construction that follows is sufficiently general to accommodate most of the mesh-based discretizations like continuous and discontinuous Galerkin finite element methods (CGFEM and DGFEM), virtual element method (VEM), finite volume method (FVM), finite difference method (FDM), and variations of these. Following the notation introduced above, we assume
\begin{equation*}
    \dim(V_h)=m
\end{equation*}
and denote a basis of $V_h$ by $\{\phi_i\}_{i=1}^m$. Therefore, 
\begin{equation*}
    V_h = \operatorname{span}\{\phi_1,\ldots,\phi_m\}.
\end{equation*}
Since $V_h$ is isomorphic to $\mathbb{R}^m$, every $u_h\in V_h$ is uniquely identified as a vector in $\mathbb{R}^m$ by the linear combination of a basis of $V_h$, namely
\begin{equation*}
    \bm{u} := [u^1,\ldots,u^m]^{\top}
\end{equation*}
such that
\begin{equation}
    u_h(x) = \sum_{j=1}^m u^j\,\phi_j(x).
    \label{eq:isomorphism}
\end{equation}

The \textbf{\textit{Gram matrix}} associated with the basis and the chosen inner product is defined by
\begin{equation}
    (\bm M)_{ij} := (\phi_i,\phi_j)_V,
    \qquad i,j=1,\ldots,m.
    \label{eq:gram-matrix}
\end{equation}
Since the basis is linearly independent, $\bm M$ is symmetric positive definite. If $u_h,v_h\in V_h$ have coefficient vectors $\bm u,\bm v\in\mathbb R^m$, respectively, then
\begin{equation}
    \label{eq:weighted-inner-product}
    (u_h,v_h)_V = \bm{u}^{\top}\bm{M}\bm{v} =: (\bm{u},\bm{v})_{\bm{M}}.
\end{equation}
The term ``Gram matrix" is used because the inner product is general. When $(\cdot,\cdot)_V$ is the $L^2(\Omega)$ inner product, $\bm M$ is the usual \textbf{\textit{mass matrix}}. In general, the Gram matrix is not diagonal even in a discrete setting, and this depends on the underlying discretization setting (CGFEM, DGFEM, FVM, FDM). In practice, we do not need exact orthogonality; discrete orthogonality suffices to approximate the exact inner product $(\cdot,\cdot)_V$. The ideal situation for the computation is using a diagonal Gram matrix. This is why mass-lumping is ideal when a numerical method leads to a non-diagonal Gram matrix.

The index partition $\mathcal P$ groups coefficient vectors, and therefore the rows of matrices whose columns are coefficient vectors, into spatially associated blocks. Let $\{a_{h,\ell}\}_{\ell=1}^n\subset V_h$ be a collection of snapshots. The snapshots may correspond to different times, parameter values, or stochastic realizations. Let $\bm a_\ell\in\mathbb R^m$ be the coefficient vector of $a_{h,\ell}$ in the basis $\{\phi_i\}_{i=1}^m$. We collect these vectors as the columns of the \textbf{\textit{snapshot matrix}}
\begin{equation}
    \bm A := [\bm a_1,\bm a_2,\ldots,\bm a_n]
    \in\mathbb R^{m\times n}.
    \label{eq:snapshot-matrix}
\end{equation}
\begin{remark}
    The interpretation of a degree of freedom depends on the spatial discretization. For example, a scalar degree of freedom may be a nodal value, a modal coefficient, or a cell average. The following analysis requires only the coefficient representation, the Gram matrix, and the index partition. It does not require associating every degree of freedom with a unique physical point.
\end{remark}

\section{Local Bases}
\label{sec:localsurrogates}
For a given index partition $\mathcal P$, let
\begin{equation*}
    m_k:=|I_k|
\end{equation*}
and define the local data matrix
\begin{equation*}
    \bm A_k := \bm A[I_k,:] \in\mathbb R^{m_k\times n}
\end{equation*}
by selecting the rows of $\bm A$ whose indices belong to $I_k$. We fix an ordering of the indices in each $I_k$ and retain the mapping between the local row indices of $\bm A_k$ and the corresponding global row indices of $\bm A$. Note that the order in which the integers in $I_k$ are listed is irrelevant as long as we keep track of the mapping.

For any collection of local matrices $\bm B_k\in\mathbb R^{m_k\times n}$, the notation
\begin{equation*}
    \vstack_{k=1,\ldots,N_p}\bm B_k
\end{equation*}
denotes the unique matrix $\bm B\in\mathbb R^{m\times n}$ obtained by placing the rows of each $\bm B_k$ at the global row indices in $I_k$. In other words, $\bm{B}[I_k,:] = \bm{B}_k$ with $k=1,\ldots,N_p$. In particular,
\begin{equation}
    \label{eq:vstack}
    \bm A = \vstack_{k=1,\ldots,N_p}\bm A_k.
\end{equation}
By construction, each index set $I_k$ is associated with the physical subdomain $\Omega_k$ introduced in Section~\ref{sec:meshdescr}. For discretizations with element-local DoFs, the rows of $\bm A_k$ correspond directly to coefficients associated with cells in $\Omega_k$. For conforming discretizations, the support of a basis function may cross a subdomain interface, but the ownership rule introduced in Section~\ref{sec:meshdescr} still assigns each algebraic index to exactly one set $I_k$. Thus, $\bm A$ should be regarded as a structured snapshot matrix: its rows are spatially associated DoFs, whereas its columns correspond to time instances, parameter values, or stochastic realizations. The graph-partitioning hierarchy exploits this spatial organization when constructing candidate local approximations.

In what follows, we first consider the ordinary singular value decomposition, which gives an exact coefficient-space Frobenius-error identity. We then consider the proper orthogonal decomposition associated with a diagonal symmetric positive-definite Gram matrix, for which we derive an exact weighted-error identity and a coefficient-space Frobenius-error bound. A central tool in the analysis is the Eckart--Young--Mirsky theorem, which we recall here for notational clarity.
\begin{theorem}[Eckart--Young--Mirsky theorem~\cite{EckartYoung1936,Mirsky1960}]
\label{thm:Eckart-Young-Mirsky}
Let $\bm A\in\mathbb R^{m\times n}$ have the compact singular value decomposition
\begin{equation*}
    \bm A
    =
    \bm U\bm\Sigma\bm V^{\top}
    =
    \sum_{i=1}^{r}
    \sigma_i\bm u_i\bm v_i^{\top},
\end{equation*}
where $r:=\rank(\bm A)$, $\bm{U}\in\mathbb R^{m\times r}$, $\bm{V}\in\mathbb R^{n\times r}$, and $\bm{\Sigma}=\operatorname{diag}(\sigma_1,\ldots,\sigma_r)$. For $r>0$, the singular values satisfy $\sigma_1\geq\cdots\geq\sigma_r>0$. When $r=0$, the compact decomposition and sums over singular values are understood to be empty. For $p\in\{0,1,\ldots,r\}$, define
\begin{equation*}
    \bm{A}^{(p)} := \sum_{i=1}^{p}\sigma_i\bm{u}_i\bm{v}_i^{\top}.
\end{equation*}
Then
\begin{equation}
\label{eq:optimalsvd}
\bm{A}^{(p)} \in
\operatorname*{argmin}_{\substack{
\bm{B} \in \mathbb{R}^{m \times n} \\
\operatorname{rank}(\bm{B}) \leq p
}}
\|\bm{A}-\bm{B}\|
\end{equation}
for every unitarily invariant norm. In particular, for $p<r$, 
\begin{equation*}
    \|\bm A-\bm A^{(p)}\|_F^2 = \sum_{i=p+1}^{r}\sigma_i^2,
\end{equation*}
and
\begin{equation*}
    \|\bm A-\bm A^{(p)}\|_2 = \sigma_{p+1}.
\end{equation*}
\end{theorem}

In \eqref{eq:optimalsvd} we use the membership symbol because a best rank-$p$ approximation might not be unique. In particular, for $1\leq p<r$, non-uniqueness can occur when a repeated singular-value block straddles the truncation index, $\sigma_p=\sigma_{p+1}$~\cite{LiStrang2020}.

In the rest of the text, we focus on the Frobenius norm because its squared error is additive over disjoint row blocks, as required for the NRMSE analysis. In contrast, the squared spectral norm of a vertically stacked residual is not generally the sum of the squared spectral norms of its blocks. For other unitarily invariant norms, truncated-SVD errors are also determined by the discarded singular values through the corresponding symmetric gauge functions, but the same additive squared-error identity need not hold. 

\subsection{Local singular value decomposition}
For each data matrix $\bm A_k$, we can compute the (local) singular value decomposition in thin form
\begin{equation}
    \label{eq:localsvd}
    \bm{A}_k=\bm{U}_k\,\bm{\Sigma}_k\,\bm{V}_k^{\top},
\end{equation}
where $\bm{\Sigma}_k=\operatorname{diag}(\bm{\sigma}_k)$ contains the $r_k:=\rank(\bm A_k)$ positive singular values in nonincreasing order. Consequently, $\bm{\Sigma}_k\in\mathbb R^{r_k\times r_k}$, $\bm{U}_k\in\mathbb R^{m_k\times r_k}$, and $\bm{V}_k\in\mathbb R^{n\times r_k}$. For $r_k=0$, the local factors are understood to be empty.
For $p_k\in\{0,1,\ldots,r_k\}$, define the rank-$p_k$ truncation of \eqref{eq:localsvd} by
\begin{equation}
    \label{eq:localsvd-trunc}
    \bm{A}_k^{(p_k)}:=
    \bm{U}_k[{:},{:}p_k]\,
    \bm{\Sigma}_k[{:}p_k,{:}p_k]\,
    \bm{V}_k^{\top}[{:}p_k,:] =\sum_{i=1}^{p_k} \bm{\sigma}_k[i]\,\bm{U}_k[:,i]\,\bm{V}_k^{\top}[i,:].
\end{equation}

For a fixed partition $\mathcal P$ and local ranks $\bm p=(p_1,\ldots,p_{N_p})$, let
\begin{equation*}
    \mathcal B_{\mathcal P,\bm p}
    :=\left\{\bm B\in\mathbb R^{m\times n}|\rank(\bm B[I_k,:])\leq p_k\ \text{for all }k\right\}.
\end{equation*}
In the following proposition, we present an error analysis for singular value decomposition. The equality derived is a prelude to the weighted POD error identities and bound developed in the next section, which will be crucial for the final algorithm described in the paper.

\begin{proposition}
Let $\bm{A}\in\mathbb{R}^{m\times n}$, $\mathcal{P}=\{I_k\}_{k=1}^{N_p}$ an index partition of $\{1,\ldots,m\}$ and $0\leq p_k\leq r_k$ the prescribed local ranks. Then
\begin{equation}
    \label{eq:svderrorfro}
    \left\|\bm{A}-\vstack_{k=1,\ldots,N_p}\bm{A}_k^{(p_k)}\right\|_F^2 = \sum_{k=1}^{N_p}\sum_{i>p_k}^{r_k}\bm{\sigma}_k[i]^2,
\end{equation}
and the reconstruction $\vstack_{k=1,\ldots,N_p}\bm{A}_k^{(p_k)}$ minimizes the Frobenius error over $\mathcal B_{\mathcal P,\bm p}$.
\end{proposition}

\begin{proof}
    From \eqref{eq:vstack} and the definition of Frobenius norm, 
    \begin{equation}
        \label{eq:vstacknorm}
        \left\|
        \vstack_{k=1,\ldots,N_p}\bm A_k
        \right\|_F^2
        =
        \|\bm A\|_F^2
        =
        \sum_{k=1}^{N_p}\|\bm A_k\|_F^2.
    \end{equation}
    The last equality holds because the sets $I_k$ of $\mathcal{P}$ form a disjoint partition of the row indices of $\bm A$. The Eckart--Young--Mirsky theorem applied to $\bm{A}_k$ gives
    \begin{equation}
        \label{eq:localerrorEYM}
        \|\bm A_k-\bm A_k^{(p_k)}\|_F^2=\sum_{i>p_k}^{r_k}\bm{\sigma}_{k}[i]^2,
    \end{equation}
    where $r_k=\rank(\bm A_k)$, \(k=1,\ldots,N_p\) and \(0\leq p_k\leq r_k\). For a fixed index partition \(\mathcal P\), the global approximation error is obtained by aggregating the errors of each index set
    \begin{align*}
        \left\|\bm{A}-\vstack_{k=1,\ldots,N_p}\bm{A}_k^{(p_k)}\right\|_F^2
        &\overset{\eqref{eq:vstack}}{=}\left\|\vstack_{k=1,\ldots,N_p}\bm{A}_k-\vstack_{k=1,\ldots,N_p}\bm{A}_k^{(p_k)}\right\|_F^2\\
        &=\left\|\vstack_{k=1,\ldots,N_p}\left(\bm{A}_k-\bm{A}_k^{(p_k)}\right)\right\|_F^2\\
        &\overset{\eqref{eq:vstacknorm}}{=}\sum_{k=1}^{N_p}\|\bm{A}_k-\bm{A}_k^{(p_k)}\|_F^2\\
        &\overset{\eqref{eq:localerrorEYM}}{=}\sum_{k=1}^{N_p}\sum_{i>p_k}^{r_k}\bm{\sigma}_k[i]^2.
    \end{align*}
For $\bm B\in\mathcal B_{\mathcal P,\bm p}$, the squared Frobenius error is the sum of independent local errors. The Eckart--Young--Mirsky theorem therefore implies the stated optimality in this block-rank-constrained class.
\end{proof}

\subsection{Local proper orthogonal decomposition}
In this section, we construct an orthonormal basis of functions with respect to $(\cdot,\cdot)_V$ or a discrete approximation of it. Here $\bm{M}$ is diagonal and positive definite, either realizing the exact Gram matrix or a mass-lumped version of it. In the latter case, the weighted norm below is the lumped discrete norm, not necessarily the norm induced by the exact Gram matrix. Let
\begin{equation}
    \bm{M}_k:=\bm{M}[I_k,I_k]
    \label{eq:local-gram-matrix}
\end{equation}
be the local Gram matrix associated with the index set $I_k$. Since $\bm M$ is symmetric positive definite, its principal submatrix $\bm M_k$ is also symmetric positive definite. The local proper orthogonal decomposition (POD) is obtained by computing the singular value decomposition \cite{KunischVolkwein2001}:
\begin{equation*}
    \widetilde{\bm{A}}_k
    :=
    \bm{M}_k^{1/2}\bm{A}_k 
    = 
    \widetilde{\bm{U}}_k\bm{\Sigma}_k\bm{V}_k^{\top}.
\end{equation*}
Here $\bm{\Sigma}_k$, $\bm{V}_k$, and $\bm{\sigma}_k$ refer to the SVD of the weighted matrix $\widetilde{\bm A}_k$, superseding their meaning in the preceding subsection. Since $\bm M_k^{1/2}$ is invertible, $\rank(\widetilde{\bm A}_k)=\rank(\bm A_k)=r_k$.
The local POD modes are defined by
\begin{equation}
    \label{eq:morthonormbasis}
    \bm{U}_k := \bm{M}_k^{-1/2}\,\widetilde{\bm{U}}_k.
\end{equation}
They are orthonormal with respect to the local Gram matrix because
\begin{equation*}
    \bm{U}_k^{\top}\bm{M}_k\bm{U}_k = \widetilde{\bm{U}}^{\top}_k\bm{M}^{-1/2}_k\bm{M}_k^{1/2}\bm{M}_k^{1/2}\bm{M}^{-1/2}_k\widetilde{\bm{U}}_k = \widetilde{\bm{U}}_k^{\top}\,\widetilde{\bm{U}}_k = \bm{I}_{r_k}.
\end{equation*}
The notation $\bm A_k^{(p_k)}$ is reused below for the weighted POD truncation, rather than the unweighted SVD truncation in Eq.~\eqref{eq:localsvd-trunc}. For $p_k\in\{0,1,\ldots,r_k\}$, with $\bm A_k^{(0)}:=\bm 0$, define the rank-$p_k$ local POD approximation
\begin{equation}
    \label{eq:localpod}
    \bm{A}_k^{(p_k)}:=
    \bm{U}_k[{:},{:}p_k]\,
    \bm{\Sigma}_k[{:}p_k,{:}p_k]\,
    \bm{V}_k^{\top}[{:}p_k,:] =\sum_{i=1}^{p_k} \bm{\sigma}_k[i]\,\bm{U}_k[:,i]\,\bm{V}_k^{\top}[i,:].
\end{equation}

\begin{proposition}
\label{prop:weighted-pod}
Let $\bm{A}\in\mathbb R^{m\times n}$, $\mathcal P=\{I_k\}_{k=1}^{N_p}$ be an index partition of $\{1,\ldots,m\}$, and $0\leq p_k\leq r_k$ be prescribed local ranks. Assume that $\bm{M}$ is diagonal and positive definite, with $\bm{M}_k$ defined in Eq.~\eqref{eq:local-gram-matrix}. Define
\begin{equation*}
    \bm{A}_{\mathcal P}:=\vstack_{k=1,\ldots,N_p}\bm{A}_k^{(p_k)},
    \qquad
    \|\bm{X}\|_{\bm{M},F}^2:=\operatorname{tr}(\bm{X}^{\top}\bm{M}\bm{X}).
\end{equation*}
Then the global weighted reconstruction error satisfies the exact identity
\begin{equation*}
    \|\bm{A}-\bm{A}_{\mathcal P}\|_{\bm{M},F}^2
    =\sum_{k=1}^{N_p}\sum_{i>p_k}^{r_k}\bm{\sigma}_k[i]^2,
\end{equation*}
and the reconstruction is optimal in the weighted Frobenius norm over the block-rank-constrained class:
\begin{equation*}
    \bm{A}_{\mathcal P}\in
    \operatorname*{argmin}_{\bm{B}\in\mathcal B_{\mathcal P,\bm{p}}}
    \|\bm{A}-\bm{B}\|_{\bm{M},F}.
\end{equation*}
Moreover, the ordinary Frobenius reconstruction error satisfies
\begin{equation}
    \label{eq:poderrorfro}
    \|\bm{A}-\bm{A}_{\mathcal P}\|_F^2
    \leq
    \sum_{k=1}^{N_p}\frac{1}{\lambda_{\min}(\bm{M}_k)}
    \sum_{i>p_k}^{r_k}\bm{\sigma}_k[i]^2.
\end{equation}
\end{proposition}

\begin{proof}
Since $\bm{M}$ is diagonal and the index sets are disjoint,
\begin{equation*}
    \|\bm{A}-\bm{B}\|_{\bm{M},F}^2
    =\sum_{k=1}^{N_p}\|\bm{M}_k^{1/2}(\bm{A}_k-\bm{B}[I_k,:])\|_F^2.
\end{equation*}
As $\bm{M}_k^{1/2}$ is invertible, the change of variables $\bm{C}_k=\bm{M}_k^{1/2}\bm{B}[I_k,:]$ preserves rank. Applying the Eckart--Young--Mirsky theorem to $\widetilde{\bm{A}}_k=\bm{M}_k^{1/2}\bm{A}_k$ therefore gives
\begin{align*}
    \min_{\bm{B}\in\mathcal B_{\mathcal P,\bm{p}}}
    \|\bm{A}-\bm{B}\|_{\bm{M},F}^2
    &=\sum_{k=1}^{N_p}\min_{\substack{\bm{C}_k\in\mathbb R^{m_k\times n}\\\rank(\bm{C}_k)\leq p_k}}
    \|\widetilde{\bm{A}}_k-\bm{C}_k\|_F^2\\
    &=\sum_{k=1}^{N_p}\sum_{i>p_k}^{r_k}\bm{\sigma}_k[i]^2.
\end{align*}
The local POD reconstructions in Eq.~\eqref{eq:localpod} attain these minima, proving the identity and weighted optimality. Finally, using the mixed inequality $\|\bm{C}\bm{X}\|_F\leq\|\bm{C}\|_2\|\bm{X}\|_F$,
\begin{align*}
    \|\bm{A}_k-\bm{A}_k^{(p_k)}\|_F^2
    &\leq\|\bm{M}_k^{-1/2}\|_2^2
    \|\bm{M}_k^{1/2}(\bm{A}_k-\bm{A}_k^{(p_k)})\|_F^2\\
    &=\frac{1}{\lambda_{\min}(\bm{M}_k)}
    \sum_{i>p_k}^{r_k}\bm{\sigma}_k[i]^2.
\end{align*}
Summing over the subdomains and using Eq.~\eqref{eq:vstacknorm} proves \eqref{eq:poderrorfro}.
\end{proof}

\begin{corollary}[Quasi-optimality in the Frobenius norm]
Under the assumptions of Proposition~\ref{prop:weighted-pod}, let
$\kappa_2(\bm{M}_k):=\lambda_{\max}(\bm{M}_k)/\lambda_{\min}(\bm{M}_k)$.
Then
\begin{equation*}
    \|\bm{A}-\bm{A}_{\mathcal P}\|_F
    \leq
    \sqrt{\max_{1\leq k\leq N_p}\kappa_2(\bm{M}_k)}
    \min_{\bm{B}\in\mathcal B_{\mathcal P,\bm{p}}}
    \|\bm{A}-\bm{B}\|_F.
\end{equation*}
\end{corollary}

\begin{proof}
For any $\bm{B}\in\mathcal B_{\mathcal P,\bm{p}}$, weighted optimality implies
\begin{align*}
    \|\bm{A}_k-\bm{A}_k^{(p_k)}\|_F^2
    &\leq\frac{1}{\lambda_{\min}(\bm{M}_k)}
    \|\bm{M}_k^{1/2}(\bm{A}_k-\bm{A}_k^{(p_k)})\|_F^2\\
    &\leq\frac{1}{\lambda_{\min}(\bm{M}_k)}\|\bm{M}_k^{1/2}(\bm{A}_k-\bm{B}[I_k,:])\|_F^2\\
    &\leq\kappa_2(\bm{M}_k)\|\bm{A}_k-\bm{B}[I_k,:]\|_F^2.
\end{align*}
Summation over $k$ and minimization over $\mathcal B_{\mathcal P,\bm{p}}$ give the result.
\end{proof}

\begin{remark}
    For a cell-centered finite-volume discretization with one scalar degree of freedom per cell and the $L^2(\Omega)$ inner product, the local mass matrix is diagonal:
    \begin{equation*}
        \bm M_k = \operatorname{diag}\left(|T_j|\right)_{j\in\mathcal C_k},
        \qquad
        \lambda_{\min}(\bm M_k) = \min_{j\in\mathcal C_k}|T_j|.
    \end{equation*}
    On a nonuniform mesh, using the local values $\lambda_{\min}(\bm M_k)$ can yield a less conservative coefficient-space bound than using the single global value $\lambda_{\min}(\bm M)$. The actual improvement also depends on how quickly the local singular values decay. A similar argument holds for the nodal discontinuous Galerkin discretization, where the mass matrix is discretely orthogonal and diagonal, and for the modal discontinuous Galerkin discretization with orthogonal polynomials where the mass matrix is $L^2$-orthogonal and diagonal.
\end{remark}

\subsection{Connection between Frobenius norm and NRMSE}
\label{subsec:fronrmseconnection}
In practical applications, we measure the coefficient-space reconstruction error using the \textbf{\textit{normalized root-mean-square error (NRMSE)}}. Let
\begin{equation*}
    \Delta := \max_{i,j}(\bm{A}[i,j]) - \min_{i,j}(\bm{A}[i,j]) > 0
\end{equation*}
be the range of the data matrix $\bm{A}$. With the reconstruction $\bm{A}_{\mathcal P}$ defined in Proposition~\ref{prop:weighted-pod}, the NRMSE can be written as
\begin{equation}
    \label{eq:nrmse-definition}
    \operatorname{NRMSE}
    \left(\bm{A},\bm{A}_{\mathcal P}\right)
    =
    \frac{1}{\Delta}\sqrt{\frac{1}{m\,n}\sum_{i,j=1}^{m,n}\left(\bm{A}[i,j]-\bm{A}_{\mathcal{P}}[i,j]\right)^2}
    =
    \frac{\left\|\bm{A}-\bm{A}_{\mathcal{P}}\right\|_F}{\Delta\sqrt{mn}},
\end{equation}
hence
\begin{equation*}
    \left\|\bm{A}-\bm{A}_{\mathcal P}\right\|_F^2 = \Delta^2\,m\,n\,\operatorname{NRMSE}\left(\bm{A},\bm{A}_{\mathcal P}\right)^2.
\end{equation*}
For a prescribed global tolerance $\varepsilon_{\NRMSE} \geq 0$, inequality \eqref{eq:poderrorfro} shows that
\begin{equation}
    \label{eq:err1}
    \sum_{k=1}^{N_p} \frac{1}{\lambda_{\min}(\bm{M}_k)}\sum_{i>p_k}^{r_k}\bm{\sigma}_k[i]^2 \leq \Delta^2\,m\,n\,\varepsilon_{\NRMSE}^2
\end{equation}
is a sufficient condition on the chosen $\mathcal{P}$ and truncation ranks $\{p_k\}_{k=1}^{N_p}$. A sufficient way to enforce this global condition is to assign nonnegative weights $w_k$ satisfying $\sum_{k=1}^{N_p}w_k=1$ and require, for every subdomain $k$,
\begin{equation*}
    \frac{1}{\lambda_{\min}(\bm{M}_k)}\sum_{i>p_k}^{r_k}\bm{\sigma}_k[i]^2\leq w_k\,\Delta^2\,m\,n\,\varepsilon_{\NRMSE}^2.
\end{equation*}
The smallest rank satisfying this local condition is
\begin{equation}
    \label{eq:local-nrmse-budget}
    p_k := \min\left\{p\in\{0,1,\ldots,r_k\}:\sum_{i>p}^{r_k}\bm{\sigma}_k[i]^2 \leq \lambda_{\min}(\bm{M}_k)\,w_k\,\Delta^2\,m\,n\,\varepsilon_{\NRMSE}^2 \right\}.
\end{equation}
Summing these local inequalities gives Eq.~\eqref{eq:err1}. Together with
Eq.~\eqref{eq:poderrorfro}, this yields
\begin{equation*}
    \operatorname{NRMSE}\left(\bm{A},\bm{A}_{\mathcal{P}}\right) \leq \varepsilon_{\NRMSE}.
\end{equation*}
There choice of the weights $w_1,w_2,\dots,w_{N_p}$ is not unique. For example, some subdomains might be more important than others, so the weights might vary accordingly. A smaller value of $w_k$ assigns a smaller squared-error budget to subdomain $k$ and therefore imposes a stricter local requirement. In this work, we propose a general approach motivated by our goal of distributing the squared-error budget uniformly per algebraic degree of freedom. In particular we use the size-proportional allocation $w_k=m_k/m$, and taking into account that $m=\sum_{k=1}^{N_p}m_k$ the local rank-selection rule becomes
\begin{equation}
    \label{eq:err3}
    p_k :=
    \min\left\{
    p\in\{0,1,\ldots,r_k\}:
    \sum_{i>p}^{r_k}\bm{\sigma}_k[i]^2
    \leq
    \lambda_{\min}(\bm{M}_k)\,\Delta^2\,m_k\,n\,\varepsilon_{\NRMSE}^2
    \right\}.
\end{equation}

If $\bm M_k=\bm I_{m_k}$ for every $k$, the local POD reduces to the local SVD; hence the singular values of the local POD coincide with the singular values of the local SVD, and Eq.~\eqref{eq:poderrorfro} holds with equality, falling back to \eqref{eq:svderrorfro}. In this case $\lambda_{\min}(\bm M_k)=1$. Therefore, the rank-selection conditions in Eqs.~\eqref{eq:local-nrmse-budget} and~\eqref{eq:err3} nevertheless remain inequalities because they require the reconstruction error to be not larger than the prescribed error budget.

\section{Greedy Graph Partitioning, Adaptivity and Data Compression}
\label{sec:compression}
The partitioned dimensionality-reduction approach presented in the previous sections can be viewed as compression algorithms. For a prescribed index partition $\mathcal{P}$, the compression ratio is
\begin{equation}
    \label{eq:crlocalpodsvd}
    \CR(\mathcal{P})=\frac{m\,n}{\displaystyle\sum_{k=1}^{N_p}(m_k+n+1)p_k}.
\end{equation}
Indeed, under this storage convention, the local surrogate associated with subdomain $k$ requires $m_k\times p_k$ entries for $\bm{U}_k$, $n\times p_k$ entries for $\bm{V}_k$, and $p_k$ entries for the singular values. By contrast, a global rank-$p$ singular value decomposition or proper orthogonal decomposition requires $(m+n+1)\,p$ entries and therefore has compression ratio
\begin{equation*}
    \CR_{\mathrm{SVD/POD}} = \frac{m\,n}{(m+n+1)p}.
\end{equation*}
Thus, the localized approach introduces extra factor entries that enable adaptive compression of different subdomains independently. These expressions count only the entries of the factor matrices and singular values. They do not include the storage required for the partition, index maps, preprocessing parameters, file headers, or numerical precision. However, in practice, this represents minimal computational overhead even for small data matrices.

The localized representation introduces additional factor entries because each subdomain has its own temporal factors, but it can improve compression when spatial localization substantially reduces the required local ranks. The remaining question is how to select the index partition $\mathcal P$.

\paragraph{Greedy hierarchical partitioning strategy}
We consider approximation spaces defined on general meshes. A standard approach for partitioning an unstructured mesh is to construct its cell-adjacency graph, in which each vertex represents a cell and two vertices are connected when the corresponding cells share an edge in 2D or a face in 3D. The graph is then divided according to prescribed balancing and edge-cut criteria.

Graph partitioning is commonly used to distribute a mesh among parallel processes, in which case the parts are constructed to have approximately equal computational loads while reducing communication across their interfaces. A candidate mesh partition can, for example, be generated with the \textsc{Metis} algorithm~\cite{karypis1998metis}; see also its Python interface~\cite{kloeckner2022pymetis}. Other graph-partitioning algorithms could be used instead. In the present work, graph partitioning is used only to generate balanced candidate spatial partitions. The resulting parts are also required to satisfy the connectivity assumptions introduced in Section~\ref{sec:meshdescr}.

\begin{figure}[h]
    \centering
    \includegraphics[width=0.75\linewidth]{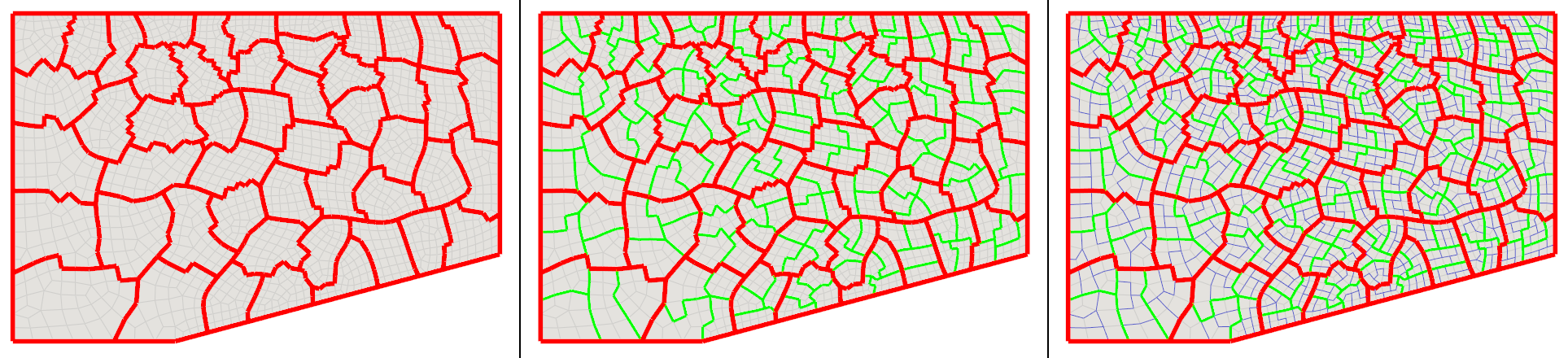}
    \caption{Example of hierarchical partition of an unstructured mesh with $N_0=50$, $N_+=4$ and $\ell=0,1,2$.}
    \label{fig:example-partitioning}
\end{figure}

The idea is to explore the space of partitions by applying the partitioning strategy hierarchically over subdomains. At the initial level, the domain $\Omega$ is partitioned into $N_0$ subdomains. At every subsequent level, a selected subdomain is tentatively divided into $N_+$ child subdomains. See Figure~\ref{fig:example-partitioning} for a visual example of the hierarchical partitioning steps. The index partition is induced by the current leaf subdomains by collecting their associated DoF indices, as described in Section~\ref{sec:meshdescr}. We use the following notation:
\begin{align*}
    \ell &:= \text{hierarchical level, with the initial partition at }\ell=0,\\
    N_0 &:= \text{number of subdomains in the initial partition},\\
    N_+ &:= \text{number of children generated by each subsequent split}.
\end{align*}

\begin{remark}
    The set of all partitions of the $N_c$ labeled cells, without imposing connectivity, has cardinality $B_{N_c}$, the $N_c$-th Bell number. The hierarchical procedure explores only the much smaller subset generated by recursive graph partitioning and satisfying the admissibility conditions introduced above. At the finest level, each subdomain consists of a single cell, which may still contain several degrees of freedom when using a high-order approximation space.
\end{remark}

\begin{remark}
Next, we present a greedy approach that starts with an initial partitioning into $N_0$ subdomains. This is justified because, for large-scale solvers for partial differential equations (parallel/distributed computing), domain partitioning is necessary to distribute computation across multiple compute nodes. However, the hierarchical partitioning can also start with no initial subdivision and recursively select partitions and local ranks.
\end{remark}

After constructing the initial partition into $N_0$ subdomains, we apply the following greedy rule recursively to each current leaf subdomain:
\begin{enumerate}[label=\Roman*.]
    \item For the current leaf subdomain $\Omega_k^{(\ell)}$, determine its
    number of rows $m_k$ and select the smallest truncation rank $p_k$
    satisfying the error condition in Eq.~\eqref{eq:err3}. Its factor-storage
    cost is
    \begin{equation*}
        S_{\mathrm{parent}} := (m_k+n+1)p_k.
    \end{equation*}

    \item Tentatively divide $\Omega_k^{(\ell)}$ into $N_+$ child subdomains
    \begin{equation*}
        \left\{\Omega_{k,i}^{(\ell+1)}\right\}_{i=1}^{N_+}.
    \end{equation*}
    For each child, determine the number of rows $m_{k,i}$ and the smallest admissible rank $p_{k,i}$ using Eq.~\eqref{eq:err3}. Define
    \begin{equation*}
        S_{\mathrm{children}}
        := \sum_{i=1}^{N_+} (m_{k,i}+n+1)p_{k,i}.
    \end{equation*}
    If $S_{\mathrm{parent}} \leq S_{\mathrm{children}}$ retain the parent as a leaf. Otherwise, accept the split and apply the same procedure recursively to its children.
\end{enumerate}

Therefore, under the size-proportional error allocation used in Eq.~\eqref{eq:err3}, replacing a parent by its children preserves the total squared-error budget assigned to that region and hence preserves the global NRMSE guarantee. Each accepted split strictly reduces the current algebraic factor storage because it replaces the parent with $N_+$ children whose combined storage is smaller. Consequently, the storage decreases monotonically along the sequence of accepted refinements. This mechanism ensures a strictly increasing algebraic compression ratio along accepted refinements.

The procedure is nevertheless greedy and does not, in general, guarantee the minimum-storage partition within the complete hierarchical search space. An immediate split may increase storage even though a deeper refinement of its children would eventually decrease it. The final partition therefore results from a sequence of locally storage-reducing decisions, rather than a globally optimal partition.

Each additional subdomain repeats temporal factors and introduces partition metadata. Larger subdomains are therefore preferable unless the reduction in local ranks is sufficient to compensate for this additional storage.

\paragraph{Loss of compression under partitioning}
For a single global domain, a rank-$p$ factorization stores $(m+n+1)p$ values.  For $N_p$ subdomains, the corresponding storage is $\sum_{k=1}^{N_p}(m_k+n+1)\,p_k$. For a balanced \textsc{Metis} partition, \(m_i \simeq m/N_p\). If all local ranks remain equal to the global rank, i.e., $p_k\simeq p$ for every $k$, then
\begin{align*}
    \sum_{k=1}^{N_p}(m_k+n+1)\,p_k
    &\simeq\sum_{k=1}^{N_p}\left(\frac{m}{N_p}+n+1\right)\,p_k\\
    &=\left(\frac{m}{N_p}+n+1\right)\,p\,N_p \\
    &=(m+N_p \, n+N_p)\,p.
\end{align*}
Hence, for $N_p>1$ and $p>0$ we have $(m+n+1)p<(m+N_p n+N_p)p$. Thus, when all local ranks remain equal to the global rank, partitioning strictly increases the factor storage relative to the global factorization. Therefore, localization is beneficial only when the reduction in the local ranks is sufficient to compensate for the repeated storage of the temporal factors.

\section{Numerical experiments}
\label{sec:numerical-results}
We assess the proposed compression method using synthetic datasets and unsteady compressible-flow simulations. The synthetic tests are designed to stress our algorithm and exhibit the mathematical properties derived in the theory, while the fluid dynamics tests consider a NACA~0012 airfoil and an inclined cube and show the algorithm's performance in practical scenarios. All meshes are generated using Gmsh~\cite{Geuzaine2009:Gmsh}. We also note that ``a posteriori NRMSE'' refers to the error between the original data and the reconstruction obtained after the compression stage, i.e., the quantity defined in equation \eqref{eq:nrmse-definition}.

\subsection{Synthetic test with solitons}
\label{subsec:synthetic-data-sets}
For this test, we consider two compactly supported traveling waves. We make this choice to perturb a prescribed spatial region, and thus a subset of spatial degrees of freedom (DOFs). Consequently, we modify only some rows of the data matrix from the zero initialization. This allows us to inject algebraic information in a controlled manner. Specifically, each traveling wave is defined by a bell-shaped function in two dimensions whose center is $(x_c,y_c)$ and the amplitude is determined by the bump function
\[
\Psi(s) =
\begin{cases}
\exp\left(\dfrac{1}{s^2 - 1}\right), & \text{if } |s| < 1, \\[6pt]
0, & \text{if } |s| \geq 1.
\end{cases}
\]
We aim to morph this function over time; hence, we define the soliton as 
\[
u(x,y,t) := A(t)\,\Psi\left(\frac{\sqrt{(x-x_c(t))^2+(y-y_c(t))^2}}{R(t)}\right),
\]
where $A(\cdot)$ is a function of time controlling the scale of the bell-shaped function, $R(\cdot)$ is a function of time controlling the radius of the circular support, and $(x_c(\cdot),y_c(\cdot))$ are functions of time controlling the center of the circular support of the bell-shaped function.

This test dataset is defined by summing two solitons. For each soliton, functions $x_c(\cdot)$, $y_c(\cdot)$, $A(\cdot)$ and $R(\cdot)$ are implemented using cubic splines. We determine the spline control points by sampling uniformly from user-defined ranges for each quantity. This approach lets us easily generate many test cases while preserving sufficient regularity in the involved functions, and it also allows us to inject information into the dataset. In Table~\ref{tab:solitons-dataset} we report the parameters necessary to reproduce the cubic spline for each quantity required to reconstruct the dataset. All splines are parametrized by the time coordinate, which spans from 0 to 1, and is divided into 500 uniform time snapshots.
\begin{table}[H]
    \centering
    \caption{Ranges, seeds, and number of points sampled from a uniform probability distribution for the parameters defining the solitons' control functions.}
    \label{tab:solitons-dataset}
    \begin{tabular}{c|ccc|ccc}
        \hline
        & \multicolumn{3}{c|}{\textbf{soliton 1}}
        & \multicolumn{3}{c}{\textbf{soliton 2}} \\
        \cline{2-7}
        parameter & range & seed & Npts
                  & range & seed & Npts \\
        \hline
        $x_c(\cdot)$ & $[1.0,9.0]$  & 45 & 10
              & $[11.0,19.0]$ & 88 & 10 \\
        $y_c(\cdot)$ & $[1.0,9.0]$  & 45 & 10
              & $[11.0,19.0]$ & 88 & 10 \\
        $R(\cdot)$   & $[1.0,1.5]$  & 50 & 20
              & $[1.4,1.5]$ & 29 & 20 \\
        $A(\cdot)$   & $[-1.0,1.0]$ & 2  & 20
              & $[-1.0,1.0]$ & 12 & 20 \\
        \hline
    \end{tabular}
\end{table}

\begin{figure}
    \centering
    \includegraphics[width=0.7\textwidth]{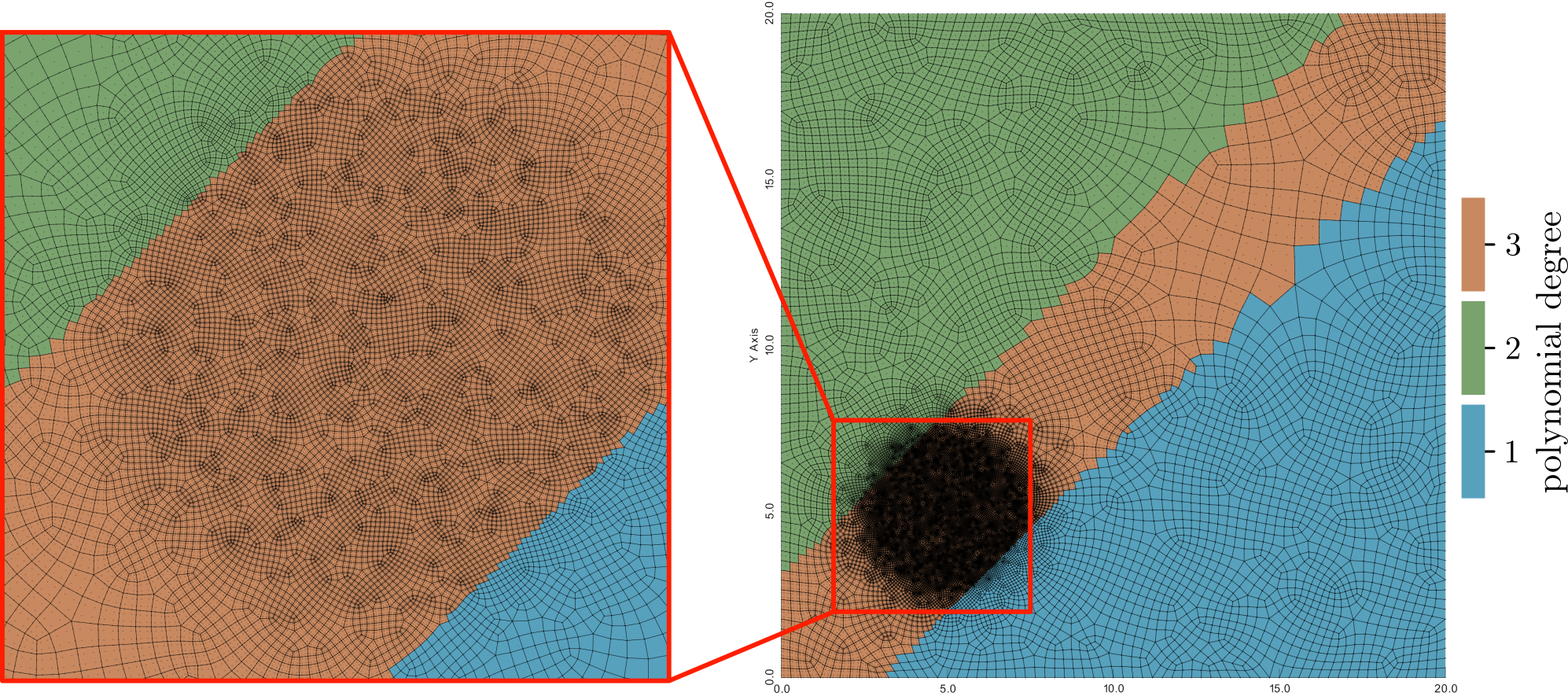}
    \caption{Mesh and polynomial degree for the solitons test.}
    \label{fig:solitons-mesh-deg}
\end{figure}

The spatial domain is a square of side 20, and the mesh contains 16,084 quadrilateral elements. The mesh is intentionally highly non-uniform. Figure~\ref{fig:solitons-mesh-deg} shows the mesh used for this test and the polynomial degree of each cell, mapped as an integer function over the domain. The corresponding LGL points in physical space are also shown in dimmed grey. The total number of points equals the number of DOFs, 206,702. Finally, the data matrix has dimensions $206,702\times500$

Since the range values for $x_c(\cdot)$ and $y_c(\cdot)$ in Table~\ref{tab:solitons-dataset} define a box where the random coordinates are sampled, the two solitons are moving in the bottom-left part and the top-right part of the square, corresponding to the refined and de-refined regions of the mesh. This configuration tends to separate the two solitons spatially; many rows in the dataset retain the initialization value, which we set to zero. 

\begin{figure}[H]
    \centering
    \begin{subfigure}[t]{0.3\textwidth}
        \centering
        \includegraphics[width=\linewidth]{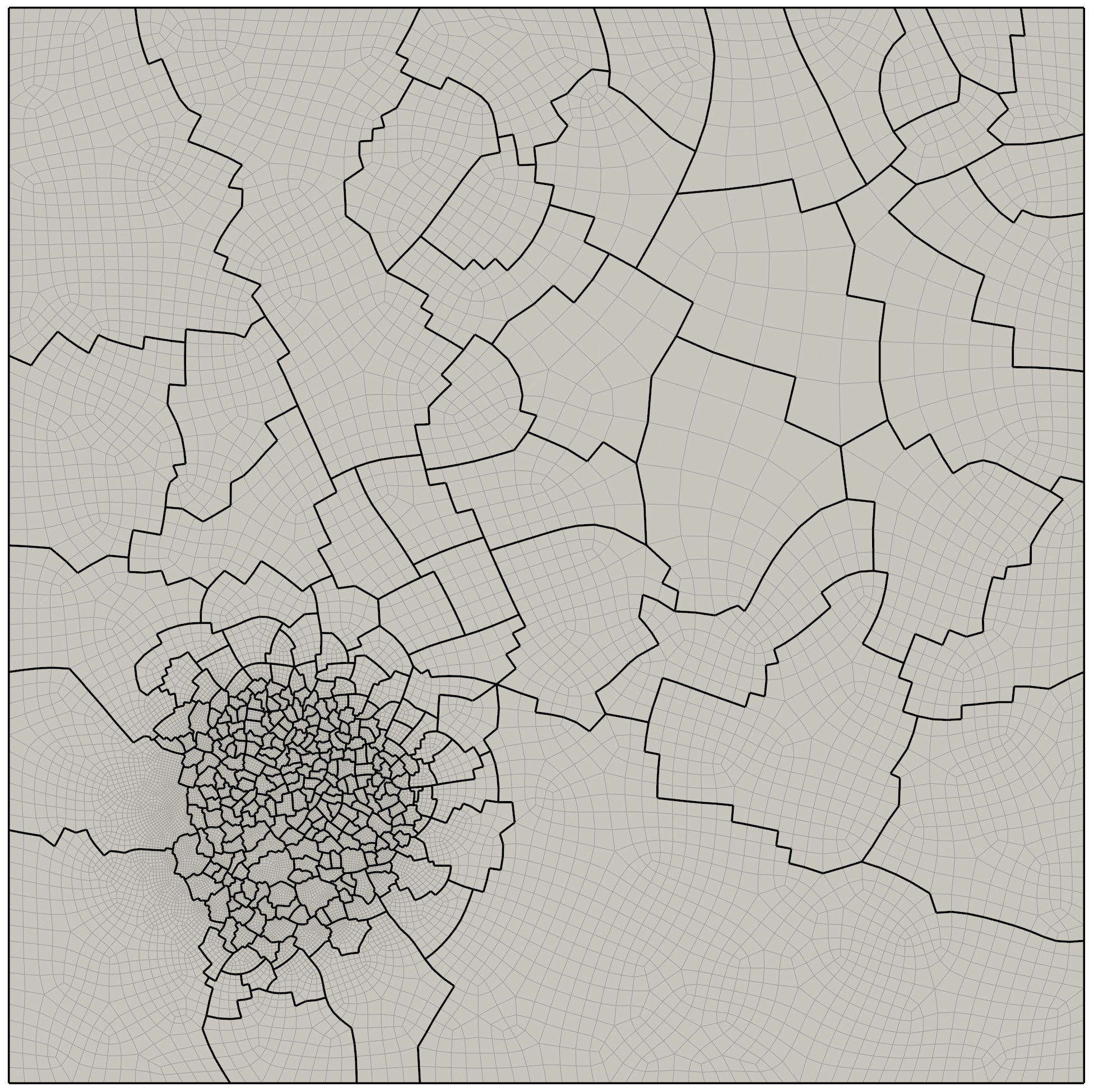}
    \end{subfigure}
    \hspace{10mm}
    \begin{subfigure}[t]{0.3\textwidth}
        \centering
        \includegraphics[width=\linewidth]{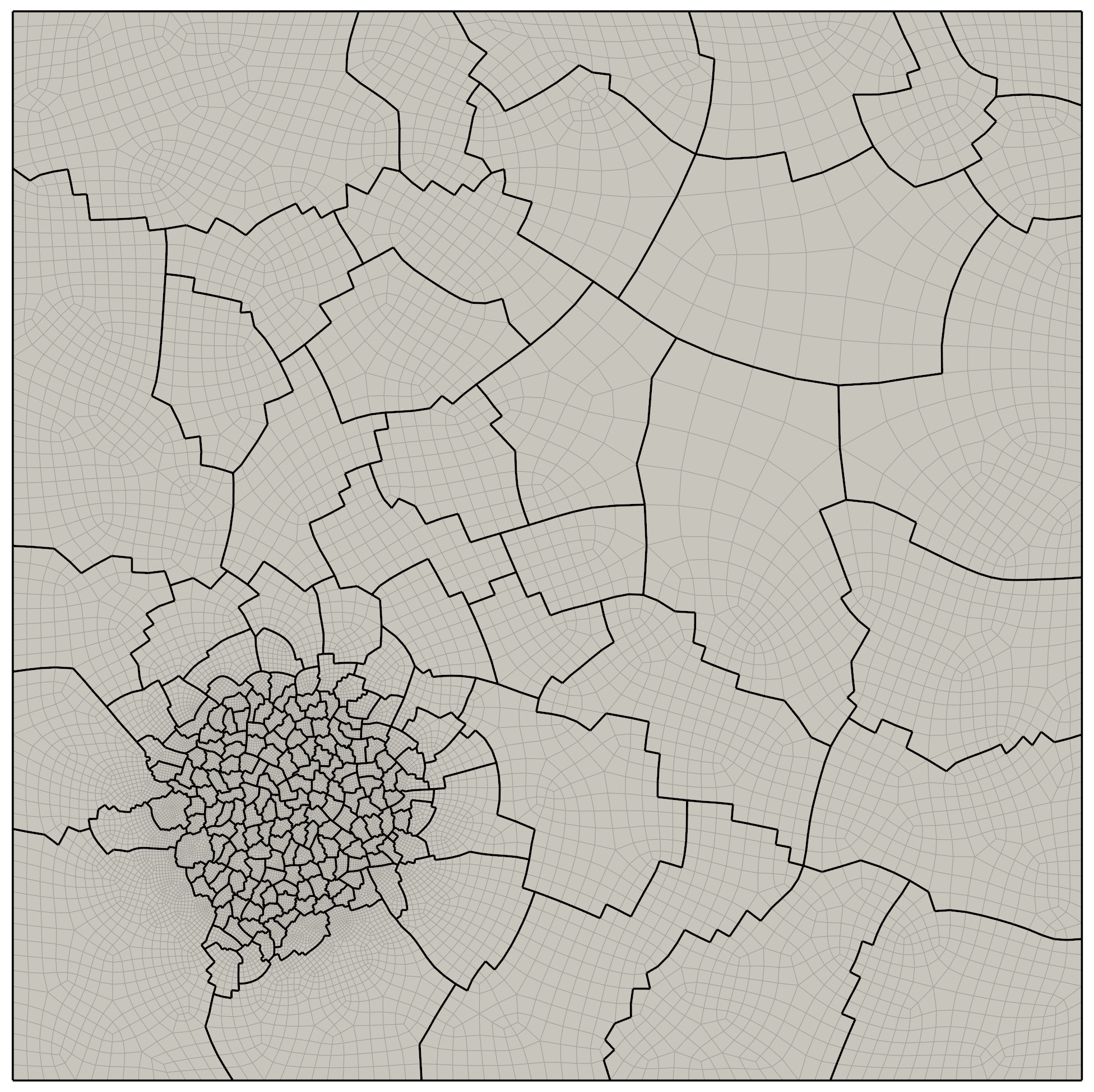}
    \end{subfigure}
    \caption{The left Figure shows the final partitioning obtained by the algorithm for the partitioning strategy 4-2 for a requested NRMSE of $10^{-4}$. The right Figure shows the final partitioning obtained by the algorithm using partitioning strategy 4-4 for a requested NRMSE of $10^{-4}$.}
    \label{fig:solitons-final-partitioning}
\end{figure}

\begin{figure}[H]
    \centering
    \begin{subfigure}[t]{0.4\textwidth}
        \centering
        \includegraphics[width=\linewidth]{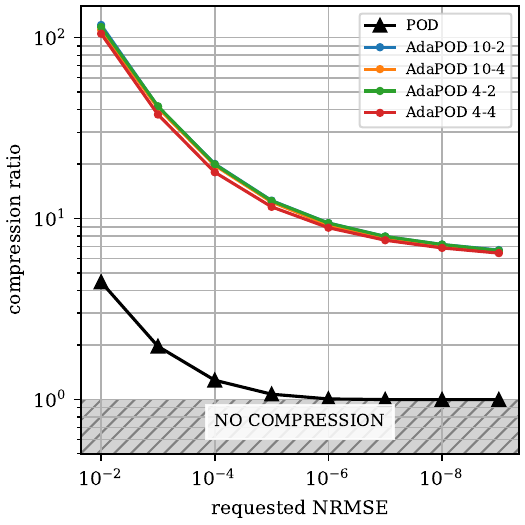}
        \caption{NRMSE vs Compression Ratio}
        \label{fig:solitons-cr}
    \end{subfigure}
    \hspace{10mm}
    \begin{subfigure}[t]{0.4\textwidth}
        \centering
        \includegraphics[width=\linewidth]{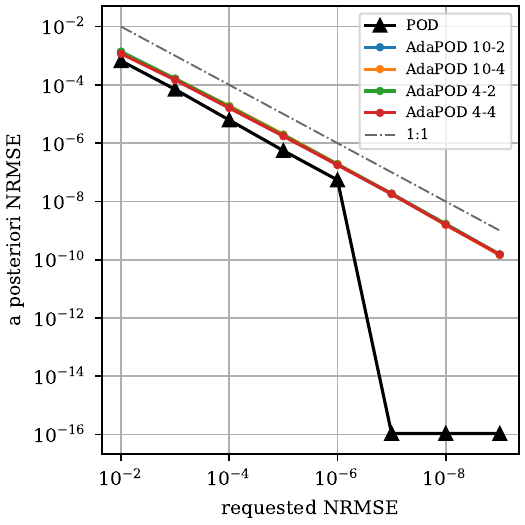}
        \caption{NRMSE requested vs NRMSE achieved}
        \label{fig:solitons-err}
    \end{subfigure}
    \caption{The figures show the algorithm's compression performance for the soliton test. In particular, Figure(a) shows the compression ratio obtained for a requested NRMSE, and Figure(b) shows the a posteriori analysis of a posteriori NRMSE  vs. requested NRMSE.}
    \label{fig:solitons-cr-err}
\end{figure}

\begin{figure}[H]
    \centering
    \includegraphics[width=0.4\textwidth]{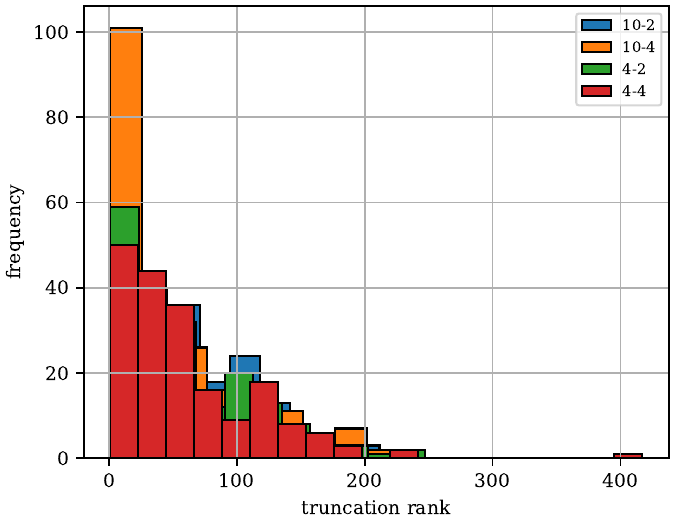}
    \caption{Frequency distribution of truncation ranks of the adaptive POD for the four partitioning strategies 10-4, 10-2, 4-2 and 4-4.}
    \label{fig:solitons-hist}
\end{figure}

In this test, we consider four partitioning trajectories: two start with 10 initial subdivisions and 2/4 subdivisions per subdomain in subsequent adaptive steps, and the other two start with 20 initial subdivisions and 2/4 subdivisions per subdomain in subsequent adaptive steps. We denote these two experiments as ``10-4", ``10-2", ``4-4" and ``4-2", respectively.
Figure~\ref{fig:solitons-cr-err} shows the compression performance for the solitons test. In particular, Figure~\ref{fig:solitons-cr} reports the compression ratio for all the partitioning strategies, and Figure~\ref{fig:solitons-err} reports the a posteriori analysis for all the partitioning strategies. The figures show that the adaptive algorithm's compression performance exceeds that of the regular POD algorithm. This is because the adaptive algorithm can localize information in space and retain only the local rank needed to satisfy the requested bound, as theory predicts. It is also interesting to observe that below a certain requested NRMSE, the POD algorithm fails to compress because the bound is too strict to allow any rank truncation; indeed, in Figure~\ref{fig:solitons-err}, for low requested NRMSE, the error drops to machine precision. In contrast, the adaptive algorithm can selectively truncate local ranks under their assigned error budgets. We also observe that the compression ratio for all partitioning strategies in Figure~\ref{fig:solitons-cr} tends toward a horizontal asymptote. The existence of a lower bound on the compression ratio does not contradict the theory. Indeed, when we set the requested NRMSE to 0, we ask the algorithm to detect the true rank of each portion of the dataset without truncation. So, if the dataset is truly low rank, we can reduce the total dimensionality without reconstruction error. Figure~\ref{fig:solitons-hist} shows the frequency of truncation ranks for the subdomains obtained at the end of the adaptive steps. As expected from theory, the adaptive algorithm tends to partition the domain into small subdomains and search for the lowest possible truncation rank, leading to clustering on the left side of the plot. Finally, Figure~\ref{fig:solitons-final-partitioning} shows the final partitioning obtained for two partitioning strategies and a requested NRMSE of $10^{-4}$. We observe that the refinement adaptivity is localized in the bottom-left part of the square, which coincides with the region with more points and variations of the signal. 


\subsection{Unsteady compressible flow tests}
\label{subsec:compressible-flow-tests}
The airfoil and cube simulations solve the compressible Navier--Stokes equations using \textsc{SSDC}, a parallel high-order flow solver~\cite{ParsaniEtAl2021SSDC}. Its spatial discretization combines discontinuous Galerkin (DG) approximations with summation-by-parts (SBP) operators and simultaneous approximation terms (SATs). The SBP property reproduces integration by parts at the discrete level, while SATs impose interelement coupling and boundary conditions weakly. This DG--SBP--SAT framework provides entropy-stable spatial discretizations on unstructured meshes. In both flow cases, compression is applied to density snapshots.

Both simulations use the explicit Bogacki--Shampine embedded Runge--Kutta pair of orders three and two~\cite{BogackiShampine1989}. The difference between the two approximations provides the error estimate used for adaptive step-size control~\cite{RanochaEtAl2022OptimizedRK}.

Dimensional quantities carry the superscript $\mathrm{dim}$; unmarked quantities are nondimensional. Length is scaled by $L_0^{\mathrm{dim}}$, chosen as the nominal airfoil chord $c^{\mathrm{dim}}$ or the cube edge $l_e^{\mathrm{dim}}$. Density and temperature are scaled by their undisturbed free-stream values, and velocity by the free-stream speed $U_0^{\mathrm{dim}}=\lvert\mathbf{u}_\infty^{\mathrm{dim}}\rvert$. Time is scaled by $L_0^{\mathrm{dim}}/U_0^{\mathrm{dim}}$. Consequently,
\begin{equation}
  \rho_\infty=T_\infty=1,
  \qquad
  \mathbf{u}_\infty=\mathbf{e}_x,
\end{equation}
where $\infty$ denotes the undisturbed free-stream state and $\mathbf{e}_x$ is the unit vector in the positive $x$ direction.

Both flow cases use a free-stream Prandtl number of $\mathrm{Pr}_\infty=0.72$. Far-field boundary conditions are imposed on the outer computational boundaries.

\subsubsection{Subsonic flow past a NACA 0012 airfoil}
\label{subsec:naca0012-case}

We consider two-dimensional unsteady flow past a NACA~0012 airfoil. The free-stream Mach number and reference Reynolds number are
\begin{equation}
  \mathrm{Ma}_\infty = 0.4,
  \qquad
  \mathrm{Re}_\infty = 5\cdot10^4.
  \label{eq:naca0012-flow-parameters}
\end{equation}

The airfoil has nominal chord $c=1$ and is rotated clockwise by $5^\circ$ about the nominal trailing edge at $(1,0)$. The mesh represents a slightly truncated trailing edge with thickness $0.002c$. The C-shaped outer boundary consists of an upstream semicircle of radius $7.3c$, centered at the midpoint of the rotated nominal chord, and horizontal upper and lower segments connected to a vertical downstream boundary at $x=6$.

No-slip adiabatic wall conditions are imposed at the top and the bottom airfoil surfaces.

The mesh contains $34{,}028$ quadrilateral elements. The solution polynomial degree is uniform throughout the domain, with $p_{\mathrm{s}}=3$. The prescribed final simulation time is $t=40$.

\begin{figure}[H]
    \centering
    \includegraphics[width=0.6\linewidth]{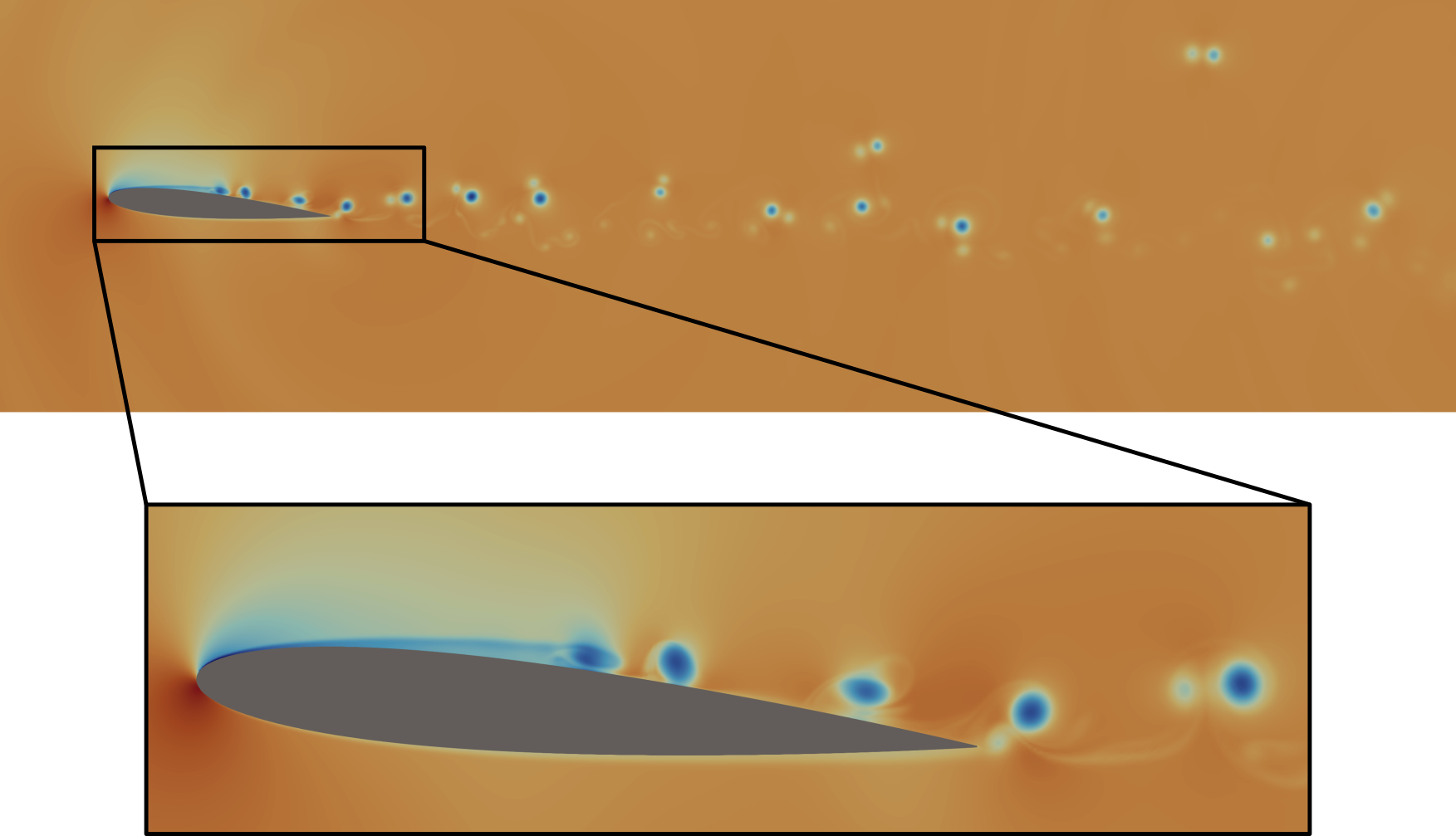}
    \caption{First snapshot of the density field sequence used to test the adaptive algorithm.}
    \label{fig:naca0012-snapshot}
\end{figure}

\begin{figure}[H]
    \centering
    \includegraphics[width=0.75\linewidth]{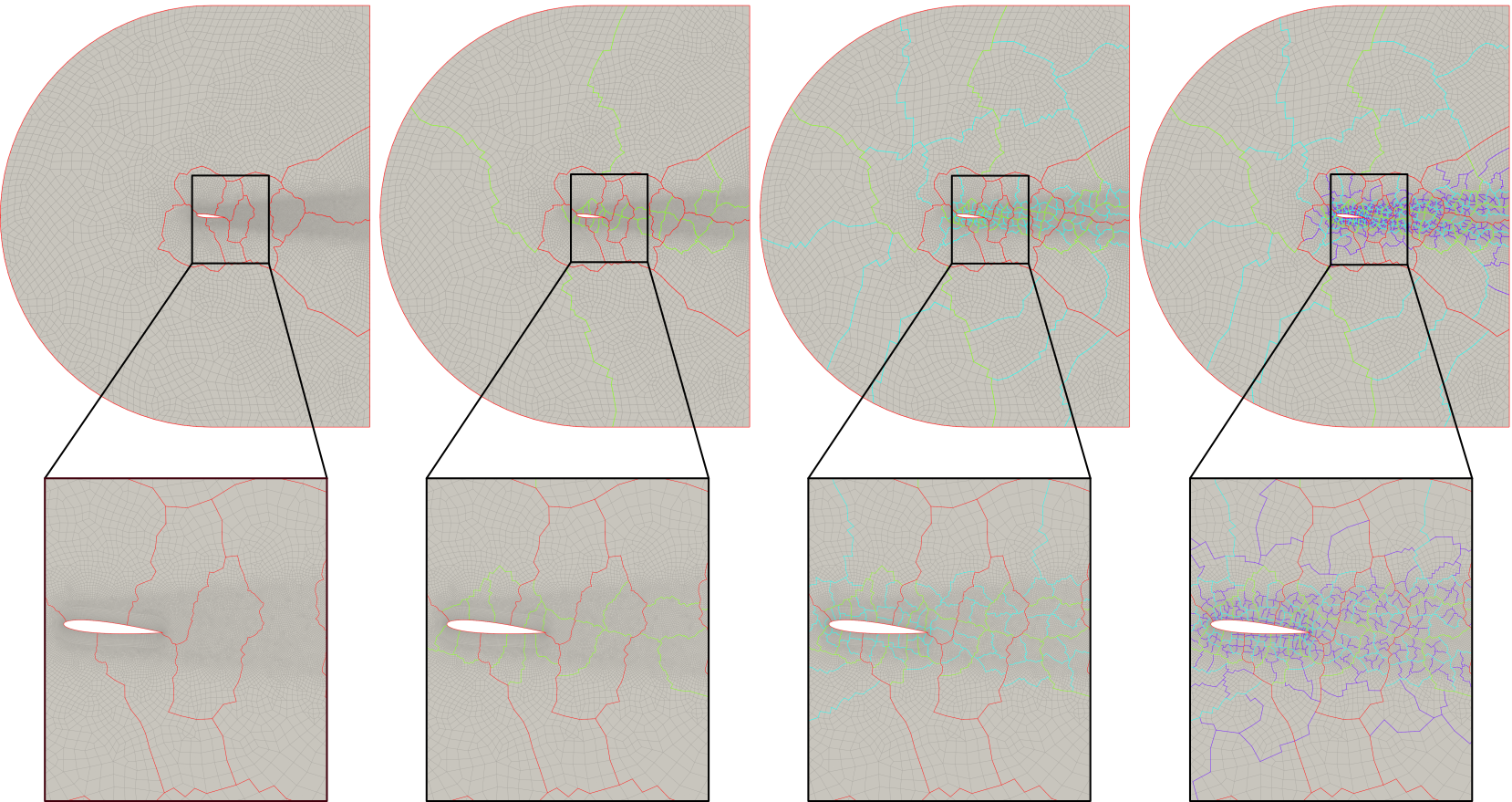}
    \caption{Four subsequent levels of adaptive partitioning for the case. The partitioning starts with 10 subdivisions, and the subsequent levels use 4 subdivisions per subdomain (see the labels ``10-4" in Figures~\ref{fig:naca0012-cr-err}) and Figure~\ref{fig:naca0012-hist}.}
    \label{fig:naca0012-part}
\end{figure}

\begin{figure}[H]
    \centering
    \begin{subfigure}[t]{0.4\textwidth}
        \centering
        \includegraphics[width=\linewidth]{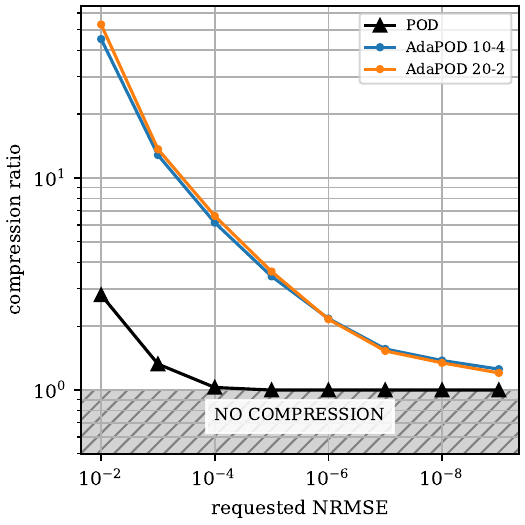}
        \caption{NRMSE vs Compression Ratio}
        \label{fig:naca0012-cr}
    \end{subfigure}
    \hspace{10mm}
    \begin{subfigure}[t]{0.4\textwidth}
        \centering
        \includegraphics[width=\linewidth]{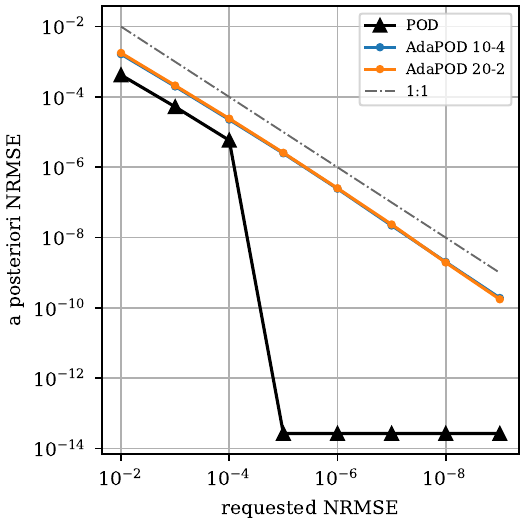}
        \caption{NRMSE requested vs NRMSE achieved}
        \label{fig:naca0012-err}
    \end{subfigure}
    \caption{The figures show the compression performance of the algorithm in the case of 500 snapshots. In particular, Figure(a) shows the compression ratio obtained for a requested NRMSE, and Figure(b) shows the a posteriori analysis of a posteriori NRMSE  vs. requested NRMSE.}
    \label{fig:naca0012-cr-err}
\end{figure}

\begin{figure}[H]
    \centering
    \includegraphics[width=0.4\textwidth]{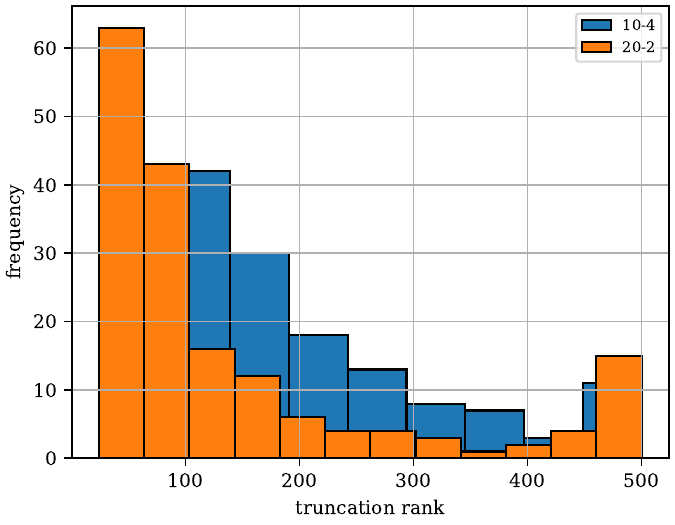}
    \caption{Frequency distribution of truncation ranks of the adaptive POD for the two partitioning strategies 10-4 and 20-2.}
    \label{fig:naca0012-hist}
\end{figure}

The airfoil data set contains $500$ density snapshots. Figure~\ref{fig:naca0012-snapshot} shows the first snapshot together with a close-up of the airfoil and its near wake. The density field exhibits localized structures near the upper surface and in the downstream wake. We start the compression when the flow is already developed to highlight the adaptive algorithm's localization/refinement effect in the wake. In this test, we consider two partitioning trajectories: one starting with 10 initial subdivisions and 4 subdivisions per subdomain in subsequent adaptive steps, and another starting with 20 initial subdivisions and 2 subdivisions per subdomain in subsequent adaptive steps. We denote these two experiments as ``10-4" and ``20-2", respectively. Figure~\ref{fig:naca0012-part} shows how the algorithm explores the refinement stages using the adaptive approach for the case ``10-4". The algorithm stops at the maximum number of refinement levels, which is 4 in this case. Figure~\ref{fig:cube-cr-err} summarizes the compression performance of our adaptive algorithm in comparison with the regular POD approach. In particular, Figure~\ref{fig:naca0012-cr} indicates very good compression properties for the requested NRMSE values. However, for very low requested NRMSE values, POD fails to compress because the required rank to satisfy the error bound equals the data matrix's maximum rank. This effect is reflected in the a posteriori analysis in Figure~\ref{fig:naca0012-err}, where the POD a posteriori error drops to zero (at machine precision). Instead, our adaptive approach shows that the a posteriori error remains very close to the requested error, both satisfying the bound and preserving compression capabilities. Finally, Figure~\ref{fig:naca0012-hist} reports the truncation ranks experienced by the local PODs computed during the adaptive stages. The clustering of rank values around the left part of the plot indicates that the adaptivity correctly identifies spatial regions containing a small amount of information; the values of ranks corresponding to the right part of the plot correspond to spatial regions rich in information, i.e., spatial regions experiencing high fluctuation in time, where the truncation rank has to be high to satisfy the error bound requirement. The difference between the frequency distributions of ``10-4" and ``20-2" is explained by the fact that ``20-2" explores more partitions than the case ``10-4" because each subdomain is divided in two, however, leading to a potentially higher number of subdomains and subsequent refinements. We treat the initial number of subdivisions and the number of subdomain subdivisions as user-defined inputs, guided by practical implementation considerations for high-performance computing scenarios. 


\subsubsection{Supersonic flow past an inclined cube}
\label{subsec:supersonic-cube-case}

We consider three-dimensional supersonic flow past a cube, following the configuration of Parsani et al.~\cite{ParsaniEtAl2026IsothermalWalls}. The free-stream Mach number and Reynolds number are
\begin{equation}
  \mathrm{Ma}_\infty=1.5,
  \qquad
  \mathrm{Re}_\infty=10^4,
  \label{eq:cube-flow-parameters}
\end{equation}
with the Reynolds number based on the cube edge.

The cube has edge length $l_e=1$ and is centered at the origin. Starting from an orientation with its edges parallel to the Cartesian axes, it is rotated by $\pi/4$ about the fixed $z$ axis and then by $\pi/4$ about the fixed $y$ axis. After these rotations, one body diagonal makes the angle
\begin{equation}
  \theta_{\mathrm{cube}}
    = \frac{1}{2}\arcsin\!\left(\frac{1}{3}\right)
    \approx 9.74^\circ
\end{equation}
with the streamwise $x$ axis, so that the incoming flow encounters an upstream vertex. The outer computational boundary is a sphere concentric with the cube, with radius $20l_e=20$.

The cube surface is stationary, impermeable, no-slip, and isothermal, with
\begin{equation}
  \mathbf{u}_w=\mathbf{0},
  \qquad
  T_w=T_\infty=1.
\end{equation}
These conditions are imposed using the complete-residual SAT (CR--SAT) formulation of Parsani et al.~\cite{ParsaniEtAl2026IsothermalWalls}.

A butterfly-type mesh contains $96{,}000$ hexahedral elements, with strong grading toward the cube. Each geometric line is divided into twenty segments with a progression factor of $1.125$. The solution polynomial degree $p_{\mathrm{s}}$ varies by region, from $p_{\mathrm{s}}=1$ near the spherical outer boundary to $p_{\mathrm{s}}=4$ near the cube.

The computation is initialized with the uniform free-stream state at $t=0$ and advanced to $t=50$. No sponge layer, artificial viscosity, shock-capturing or shock-tracking mechanism, or other shock-specific treatment is used.

\begin{figure}[H]
    \centering
    \includegraphics[width=0.6\linewidth]{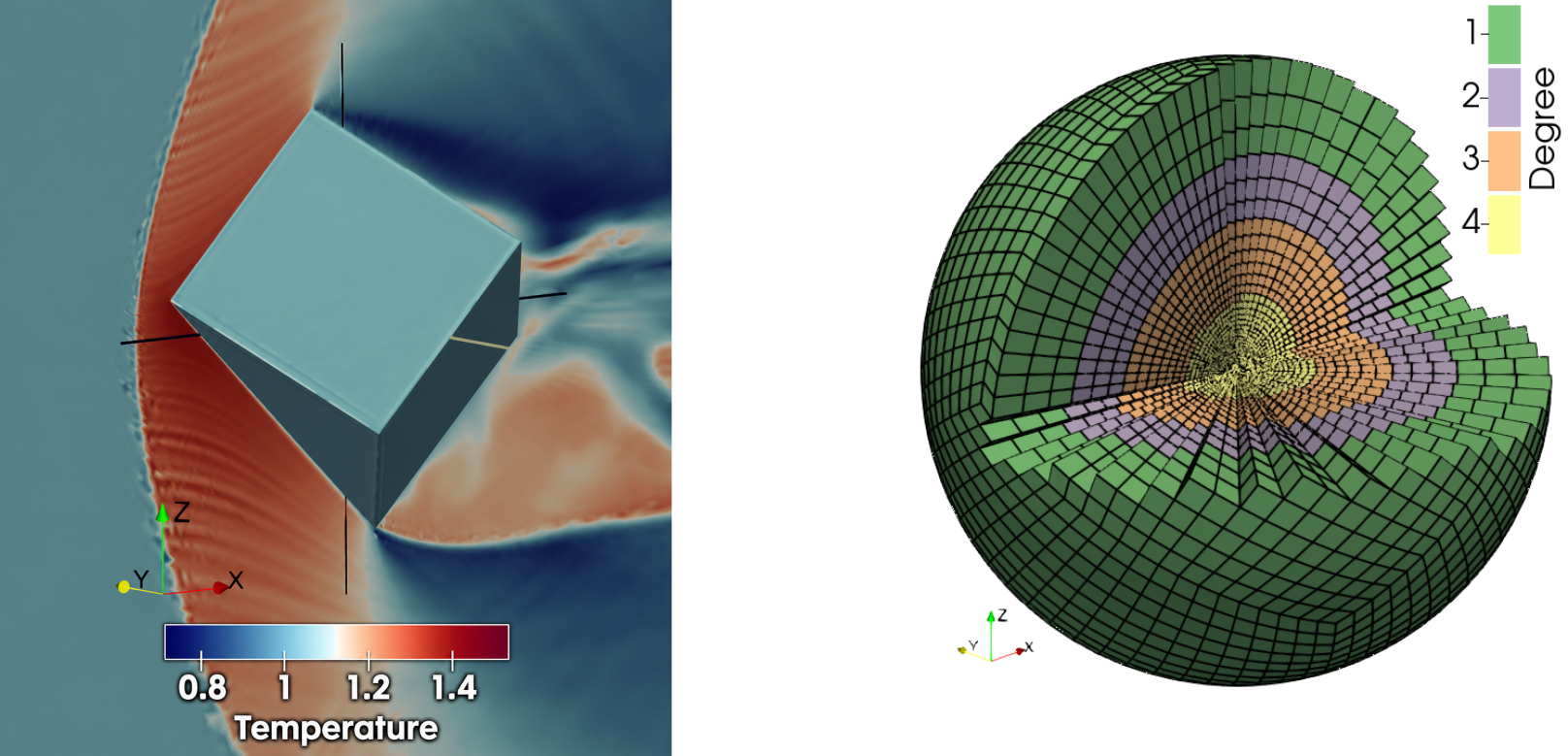}
\caption{Supersonic flow past the inclined cube: temperature in the $x$--$z$ plane at $t=50$ (left) and spatial distribution of the solution polynomial degree (right).}
    \label{fig:cube-sol-deg}
\end{figure}

Figure~\ref{fig:cube-sol-deg} shows the temperature field near the cube in the $x$--$z$ plane at $t=50$ and the spatial distribution of the solution polynomial degree. The temperature field displays a detached bow shock and a heated windward region.

\begin{figure}[H]
    \centering
    \includegraphics[width=0.75\linewidth]{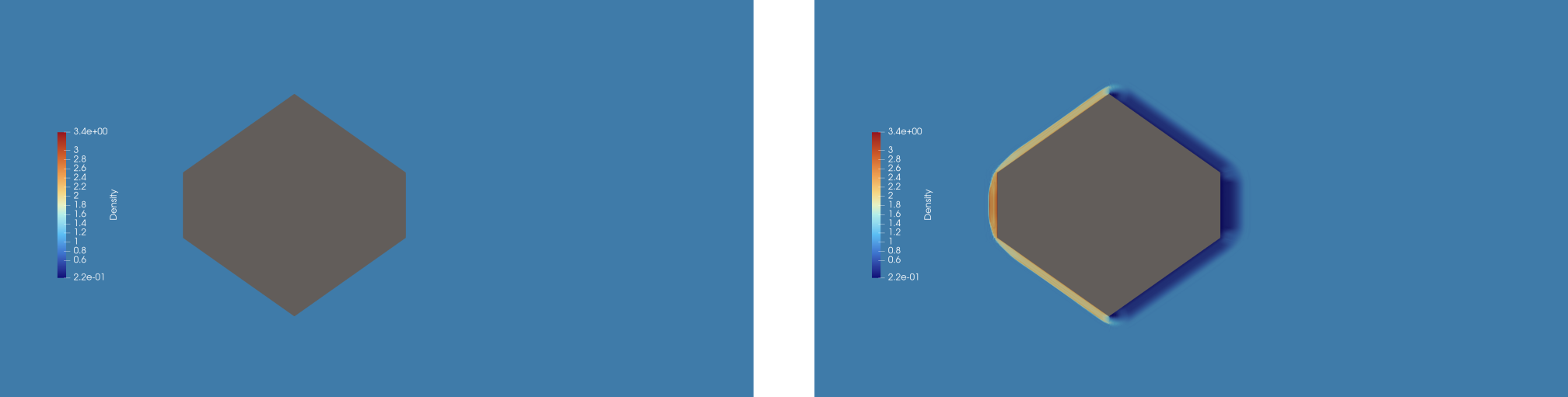}
\caption{Supersonic flow past the inclined cube, dataset 0: on the left, the density field at initial time with a uniform spatial value; on the right, the density field after 100 timesteps.}
    \label{fig:cube-0-99}
\end{figure}

\begin{figure}[H]
    \centering
    \includegraphics[width=0.75\linewidth]{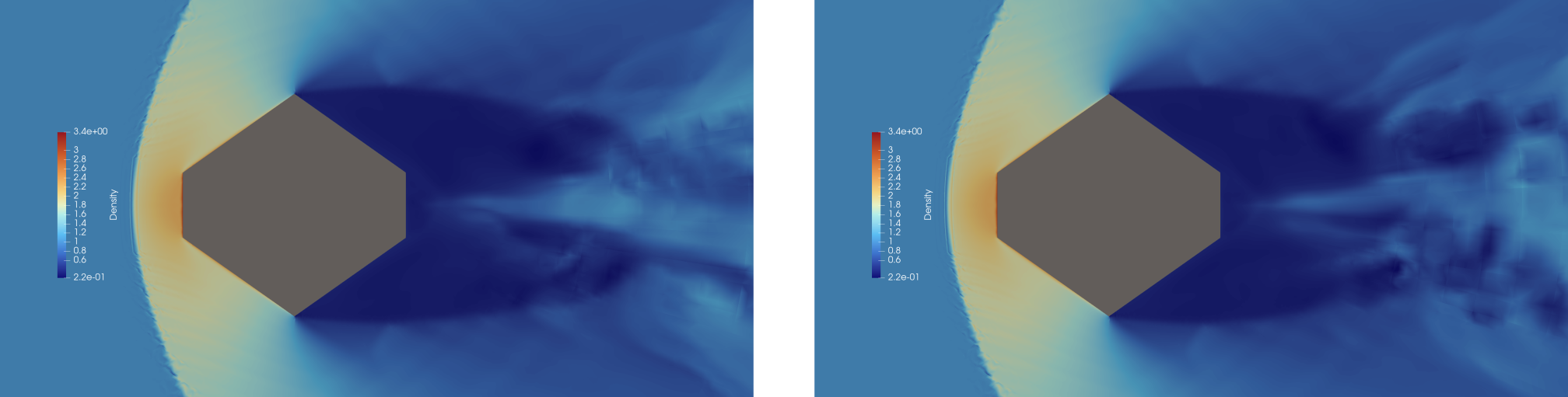}
\caption{Supersonic flow past the inclined cube, dataset 66,000: on the left, the density field at step 66,000; on the right, the density field after 100 timesteps.}
    \label{fig:cube-66000-99}
\end{figure}

\begin{figure}[H]
    \centering
    \begin{subfigure}[t]{0.4\textwidth}
        \centering
        \includegraphics[width=\linewidth]{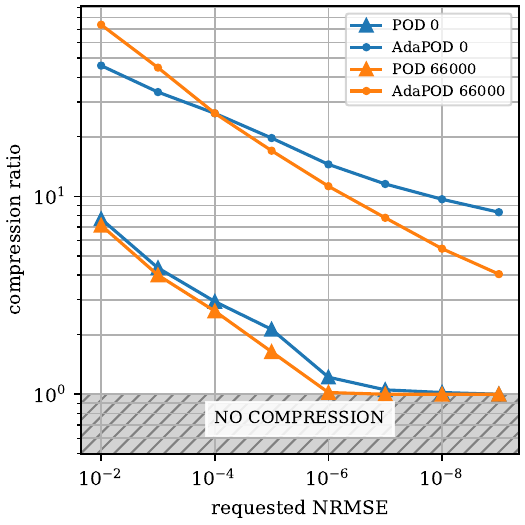}
        \caption{NRMSE vs Compression Ratio}
        \label{fig:cube-cr}
    \end{subfigure}
    \hspace{10mm}
    \begin{subfigure}[t]{0.4\textwidth}
        \centering
        \includegraphics[width=\linewidth]{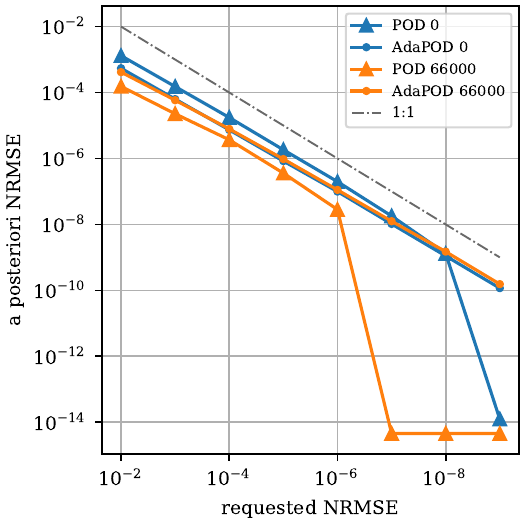}
        \caption{NRMSE requested vs NRMSE achieved}
        \label{fig:cube-err}
    \end{subfigure}
    \caption{The figures show the algorithm's compression performance for 100 snapshots; the numbers 0 and 66000 indicate the starting timestep used to create the 100-snapshot dataset. In particular, Figure(a) shows the compression ratio obtained for a requested NRMSE, and Figure(b) shows the a posteriori analysis of a posteriori NRMSE  vs. requested NRMSE.}
    \label{fig:cube-cr-err}
\end{figure}

\begin{figure}[H]
    \centering
    \includegraphics[width=0.4\textwidth]{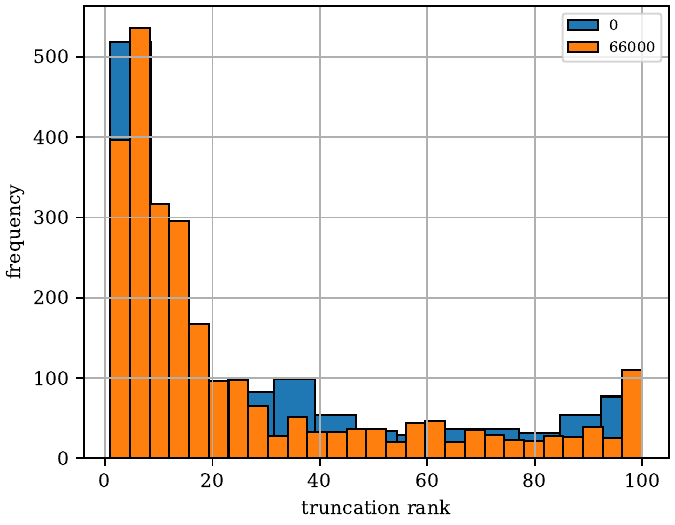}
    \caption{Frequency distribution of truncation ranks of the adaptive POD for the two tests: the label 0 in the legend indicates the dataset of 100 snapshots starting from time 0; the label 66,000 indicates the dataset of 100 snapshots starting from step 66000.}
    \label{fig:cube-hist}
\end{figure}

For the compression tests, we consider the density field analogously to the two-dimensional airfoil test. We consider two datasets: the first dataset contains the first 100 timesteps of the simulation. Figure~\ref{fig:cube-0-99} shows a section of the plane $z=0$ near the cube for the first and last timesteps. We denote such a dataset with ``0" in the legend of the following result figures. Similarly, we construct a second dataset of 100 timesteps when the flow is already fully developed. Figure~\ref{fig:cube-66000-99} shows a section of the plane $z=0$ near the cube at steps 66000 and 669,99. We denote such a dataset with ``66,000" in the legend of the following result figures.
Figure~\ref{fig:cube-cr}, Figure~\ref{fig:cube-err} and Figure~\ref{fig:cube-hist} show the compression performance for the two datasets. The comments for this case match those for the NACA~0012 test. The adaptive algorithm shows strong compression performance. Notably, the two datasets considered here differ markedly in their spatial and temporal features, which is why we selected them for testing. However, the compression performance remains consistent, experimentally proving the efficacy of our adaptive approach.

\section{Conclusions}
\label{sec:conclusions}
We have developed a spatially adaptive local proper orthogonal decomposition method for error-bounded lossy compression of snapshot data in finite-dimensional function spaces. The method combines graph-based spatial partitioning with independent low-rank approximations on the resulting subdomains. Unlike global POD, the proposed representation allows the local approximation ranks to vary across the computational domain, thereby exploiting spatial variations in the complexity of the snapshot data.

The mathematical analysis establishes an exact decomposition of the global squared weighted reconstruction error into local contributions for partitions that are orthogonal in the discrete inner product. For a fixed partition and prescribed local ranks, the corresponding local POD approximations minimize the weighted Frobenius error within the associated block-separable approximation class and are quasi-optimal in the ordinary Frobenius norm. These properties enable the allocation of a global error budget among subdomains and the selection of the smallest admissible local ranks. A sufficient bound also guarantees that the coefficient-space NRMSE does not exceed a prescribed tolerance. The resulting hierarchical algorithm accepts a refinement only when it reduces the algebraic storage required by the low-rank factors while preserving the global error guarantee.

The numerical experiments consider synthetic data and snapshot sequences from two- and three-dimensional compressible-flow simulations. They demonstrate the ability of spatially adaptive partitioning to identify regions with different approximation requirements and to allocate local ranks accordingly. The results illustrate the potential compression advantages over global POD and fixed spatial partitions, particularly for data with spatially localized or heterogeneous dynamics. The reconstruction errors remain within the prescribed tolerances, consistently with the theoretical error bounds.

The method has two important limitations. First, the partitioning strategy is greedy and does not guarantee a globally minimum-storage representation. Second, the compression ratio is based on the algebraic storage of the low-rank factors and excludes partition metadata and encoding overhead. Future work will investigate more general error-budget allocations, alternative partition search strategies, integrating local representations into distributed scientific data compression workflows, and applications to reduced-order modeling.

\section*{Declarations}
\subsection*{CRediT authorship contribution statement}
\textbf{Roberto Nuca:} Conceptualization, Formal analysis, Investigation, Methodology, Coding, Visualization, Data Generation, Data Curation, Writing - original draft, Writing - review \& editing.\\
\textbf{Matteo Parsani:} Investigation, Methodology, Supervision, Visualization, Writing - original draft, Writing - review \& editing, Funding acquisition.\\
\textbf{Lisandro Dalcin:} Coding, Data Generation.

\subsection*{Data Availability}
The datasets generated and analyzed during the current study are available from the corresponding author upon reasonable request. The synthetic datasets can be reproduced using the procedures and parameters described in the manuscript. The numerical simulation datasets, including the snapshot sequences for the two-dimensional NACA 0012 airfoil and the three-dimensional inclined-cube configurations, are available upon reasonable request.

\subsection*{Acknowledgments}
The work described in this paper was supported by King Abdullah University of Science and Technology through the grant number BAS/1/1663-01-01. The authors are also thankful for the computing resources of the Supercomputing Laboratory at King Abdullah University of Science and Technology. Roberto Nuca is a member of the GNCS/INdAM research group.

\subsection*{Conflict of interest}
The authors declare no conflict of interest.

\subsection*{AI usage statement}
The research presented here was conducted without using AI tools. The authors developed the original ideas, mathematical arguments, and code used to implement the algorithm and conduct the tests without AI assistance. We used AI-assisted tools solely to improve the manuscript's grammar and text formatting.

\bibliographystyle{plain}
\bibliography{references}
\end{document}